\documentclass[leqno]{degm}
\usepackage[colorlinks]{hyperref}

\contact[divizio@picard.ups-tlse.fr]
{Lucia Di Vizio, Laboratoire Emile Picard, Universit\'e Paul Sabatier, U.F.R. M.I.G.,
118, route de Narbonne, 31062 Toulouse CEDEX4, France}

\title{Introduction to $p$-adic $q$-difference equations\\
{\large (weak Frobenius structure and transfer theorems)}}

\author{Lucia Di Vizio}

\makeatletter
\renewcommand{\thesection}{\@Roman\c@section}
\renewcommand{\thesubsection}{\@arabic\c@subsection}
\makeatother

\theorembodyfont{\sl}
\newtheorem{thm}{Theorem}[subsection]
\newtheorem{lemma}[thm]{Lemma}
\newtheorem{prop}[thm]{Proposition}
\newtheorem{cor}[thm]{Corollary}
\newtheorem{defn}[thm]{Definition}

\theorembodyfont{\rm}
\newtheorem{rmk}[thm]{Remark}

\newtheorem{parag}[thm]{\hskip -4pt}
\newtheorem{notation}[thm]{Notation}

\numberwithin{equation}{thm}

\def\ds{\displaystyle}

\def\dq{d_q}

\def\sgq{\sigma_q}
\def\nq{[n]_q}

\def\C{{\mathbb C}}
\def\Q{{\mathbb Q}}
\def\Z{{\mathbb Z}}
\def\R{{\mathbb R}}
\def\N{{\mathbb N}}

\def\l{\left}
\def\r{\right}
\def\[[{\l[\!\l[}
\def\]]{\r]\!\r]}
\def\p{\prime}

\def\a{\alpha}
\def\be{\beta}
\def\ga{\gamma}
\def\la{\lambda}
\def\veps{\varepsilon}

\def\crash#1{}

\begin{document}

\maketitle

\begin{abstract}
Inspired by the theory of $p$-adic differential equations, this paper introduces
an analogous theory for $q$-difference equations over a local field, when $|q|=1$.
We define some basic concepts, for instance the \emph{generic radius of convergence},
introduce technical tools, such as
a \emph{twisted Taylor formula} for analytic functions, and prove some fundamental
statements, such as an \emph{effective bound theorem}, the \emph{existence of a weak Frobenius structure}
and a \emph{transfer theorem} in regular singular disks.
\end{abstract}


\renewcommand{\thefootnote}{ }
\footnotetext{The author acknowledges the Institute for the
Advanced Study for hospitality and the National Science Foundation for partial support
during the later stages of preparation of this paper.}

\tableofcontents


\section*{Introduction}


Since the late 1940's, $q$-difference equations have been almost forgotten.
In the last ten years the
field has recovered its original vitality and the theory has witnessed substantial advances.
Authors have also considered these functional equations both
from an arithmetical and a $p$-adic point of view (cf. for instance
\cite{BB}, \cite{AnnII} or \cite{DV}).
This paper seeks to fill a gap in the literature, offering a systematic
introduction to $p$-adic $q$-difference equations \emph{when $|q|=1$}.

\par
Our motivation for this work was Sauloy's result
on $q$-deformations of the local monodromy of fuchsian differential
systems over ${\mathbb P}^1_\C$
(cf. \cite{Sfourier} and, for a survey on the topic, see \cite{Aqdef}, \cite{DRSZ}).
Sauloy considers a fuchsian differential system
$$
x\frac{dY}{dx}=G(x)Y(x)\ ,
\leqno{({\mathcal S})}
$$
such that $G(x)\in M_{\mu\time\mu}(\C(x))$.
More precisely he supposes that the matrix $G(x)$ has no poles at $0$ and $\infty$,
that the difference of any two eigenvalues of $G(0)$ (resp. $G(\infty)$) is not a non-zero integer
and that all the poles $x_1,\dots,x_s$ of $G(x)$ are simple.
Under these assumptions one can choose $q\in\C$, with $|q|\neq 1$, such that
$x_iq^\R\cap x_jq^\R=\emptyset$ for every $i\neq j$ and construct for any $0<\veps<<1$
a $q^\veps$-difference system of the form
$$
Y(q^\veps x)=\Big[{\mathbb I}_\mu+(q^\veps-1)G_\veps(x)\Big]Y(x)\ ,
\leqno{({\mathcal S})_{q^\veps}}
$$
such that the matrix $G_\veps(x)\in M_{\mu\time\mu}(\C(x))$ converges uniformly to $G(x)$ over
$\Omega=\C^\ast\smallsetminus \l(\cup_i x_iq^\R\cup q^\R\r)$. One can suppose, for instance by taking $G_\veps(x)=G(x)$ for all $\veps$,
that there exists a matrix $Q_\veps\in Gl_\mu(\C)$ such that $Q=\lim_{\veps\rightarrow\infty}Q_\veps$
and both $Q_\veps G_\veps(0)Q_\veps^{-1}$ and $QG(0)Q^{-1}$ are upper triangular matrices,
and that an analogous hypothesis is verified at $\infty$.
Then, for $\veps$ sufficiently close to $0$, Sauloy constructs two canonical solutions $Y_{\veps,0}(x)$
and $Y_{\veps,\infty}(x)$ of $({\mathcal S})_q$,
respectively in a neighborhood of $0$ and $\infty$: they both turn out to be meromorphic on $\C^\ast$.
Therefore the Birkhoff matrix $P_\veps(x)=Y_{\veps,\infty}(x)^{-1}Y_{\veps,0}(x)$
is also meromorphic over $\C^\ast$ and moreover it is elliptic: $P_\veps(x)=P_\veps(q^\veps x)$.

\begin{thm}{\rm \cite[\S 4]{Sfourier}}
The matrix $P_\veps(x)$ tends to a locally constant matrix $P(x)$ over $\Omega$, when $\veps\rightarrow 0$.
Let $P^\p$ and $P^{\p\p}$ be the values of $P(x)$ over the two connected components of $\Omega$
whose closure contains the pole $x_i$ of $G(x)$.
Then the local monodromy of $({\mathcal S})$ around $x_i$ is generated by $P^\p{P^{\p\p}}^{-1}$.
\end{thm}

Christol and Mebkhout have developed a theory of \emph{$p$-adic monodromy},
answering questions first raised by Robba: due to the lack of analytic continuation
it is much more complicated than the complex theory, therefore we think it would be interesting to study
the properties of $q$-deformations in the $p$-adic setting.
\par
There is a fundamental difference between complex and $p$-adic $q$-deformations:
while in $\C$ one can let $q\rightarrow 1$ avoiding the unitary circle,
this is not possible in the $p$-adic world. In other words, to study $p$-adic $q$-deformations of
differential equations one has to deal with the case $|q|=1$.
\par
We are emphasizing the ``$|q|=1$'', since the literature is almost
entirely
devoted to the case $|q|\neq 1$, to avoid the small divisor problem.
In fact, consider a $q$-difference equation
$$
y(q^\mu x)+a_{\mu-1}(x)y(q^{\mu-1}x)+\dots+a_0(x)y(x)=0\ ,
\leqno{({\mathcal E})_q}$$
where $q$ is an element of a normed field, archimedean or not, such that
$|q|\neq 1$, and the $a_i(x)$'s are rational functions in $K(x)$.
Then any convergent solution $\sum_{n\geq 0} a_nx^n$ of $({\mathcal E})_q$ is the
expansion at $0$ of a meromorphic function over ${\mathbb A}^1_K$, since the equation itself
allows for a meromorphic continuation of the solution.
Moreover, the Adams lemma states that even the most irregular $q$-difference equations
always have at least one solution whose uniform part is analytic at zero
(cf. \cite{AD} and \cite[1.2.6]{Spoligono}).
On the other hand, when $|q|=1$, one has to deal with the problem of estimating
terms of the form $1-q^n$ (the so called small divisors problem), which can make the prediction of
the existence of a convergent solution of  $({\mathcal E})_q$ very difficult.
\par
We point out that we are not distinguishing between the archimedean and the ultrametric
case: while the differential equation theory
gives rise to two substantially different theories in the $p$-adic and the complex framework,
the dichotomy in the $q$-difference world is  given by the cases
$|q|=1$ and $|q|\neq 1$.
\par
Nowadays, the complex theory of $q$-difference equations for $|q|\neq 1$ has reached
a ``degree of completeness'' comparable to differential equation theory,
as Birkhoff and Guenther hoped \cite{Birk}. Moreover it seems that those results should also be
true in the $p$-adic case with very similar proofs
(cf. for instance  \cite{Beindex} versus \cite{BB} and \cite[I, 2.2.4]{Stesi}).
In the meantime very few pages are devoted to the study of the case $|q|=1$,
which, apart from the small divisor problem, is characterized by essentially two difficulties:
\begin{trivlist}
\item 1)
The meromorphic continuation of solutions does not work any longer, so one needs
a good notion of solution at a point $\xi\neq 0$.  This is a key-point of the $p$-adic approach, since
we cannot imagine a $p$-adic theory of $q$-difference equations without an analogue of the
notion of Dwork-Robba generic point and generic radius of convergence.
\item 2)
Sauloy constructs his canonical solutions using the classical Theta function
$$
\Theta(x)=\sum_{n\in\Z}q^{-\frac{n(n-1)}{2}}x^n\ .
$$
Of course it does not converge for $|q|=1$.
Moreover all the infinite products, which play
such an important role in the theory of $q$-series (cf. for instance \cite{GR})
and would be very useful to write meromorphic solutions, do not converge either.
\end{trivlist}
The lesson of the results contained in this paper is that $p$-adic $q$-difference equations with $|q|=1$
present the same pathologies as $p$-adic differential equations, namely the
uniform part of solutions at a regular singular points can be divergent,
according to the type of the exponents.
In fact, the $q$-difference theory is precisely a $q$-deformation of the differential situation.
Since $p$-adic differential equation theory has geometrical implications, this allows
one to imagine that $p$-adic $q$-difference equations should have some
``non commutative $p$-adic geometric implications''
\cite{Villetaneuse}.

$$
\ast\ \ast\ \ast
$$

Concerning the content of the present work, we have chosen to introduce
the notions of $q$-difference algebra as required by the paper.
Anyway a systematic presentation with highly compatible notation can be found in \cite[\S 1]{DV}.
The paper is organized as follows:

\begin{description}

\item[] Chapter \ref{sec:introduction} is an introduction to basic tools.
First of all we study the properties of the $q$-difference algebra of analytic functions over a disk
and the properties of $q$-difference operators with respect to the so called Gauss norms.
We also state a result, proved
in Appendix \ref{sec:taylor}, about the existence of a $q$-expansion of analytic functions,
which we use to establish a good notion of solutions of a $q$-difference
equation at a point $\xi\neq 0,\infty$. In particular this allows for the definition
of generic radius of convergence at a Dwork-Robba generic point.

\item In Chapter \ref{sec:ee} we prove an effective bound theorem in the wake of the Dwork-Robba theorem,
from which we deduce a transfer theorem in ordinary $q$-orbits and a corollary on the
$q$-deformations of $p$-adic differential equations.

\item Following Christol \cite{Ctransfer}, in Chapter \ref{sec:frobenius}
we construct a weak Frobenius structure for $q$-difference
systems having a regular singularity at $0$.

\item In Chapter \ref{sec:transfer} we prove a $q$-analogue of the Christol-Andr\'e-Baldassarri-Chiarel\-lotto
transfer theorem in a regular singular disk, relying on the result of Chapter \ref{sec:frobenius}.

\item The Appendix is divided in three independent parts.

Appendix \ref{sec:taylor} contains a proof the twisted Taylor formula stated in
\S\ref{subsec:anfunction} for analytic functions over a disk
and a generalization to analytic functions over non connected analytic domains.
The proof is completely elementary and the techniques used in it do not play any role in the paper.

In Appendix \ref{sec:regsing} we quickly recall some basic facts about
regular singular $q$-dif\-fe\-rence systems that we use in Chapters  \ref{sec:frobenius}
and \ref{sec:transfer}.

Finally in Appendix \ref{sec:type} we have grouped some technical estimates
of the $q$-type, which are used in Chapter \ref{sec:transfer}.

\end{description}


\section{Basic definitions and properties of $p$-adic $q$-difference systems}
\label{sec:introduction}

Until the end of the paper $K$ will be an algebraically closed field of characteristic zero,
complete with respect to a non-archimedean norm $|~|$, inducing a $p$-adic norm
over $\Q\hookrightarrow K$. We fix the normalization of $|~|$ by setting $|p|=p^{-1}$.
Moreover, we fix an element $q\in K$, such that
\begin{enumerate}
\setlength{\itemsep}{-2pt}
\item $|q|=1$;
\item $q$ is not a root of unity;
\item the image of $q$ in the residue field of $K$ generates a
finite cyclic group (i.e. it is algebraic over ${\mathbb F}_p$).
\end{enumerate}

\noindent
Let us consider the ring

$$
{\mathcal A}_D=\l\{\sum_{n\geq 0}a_n(x-\xi)^n: a_n\in K,
\liminf_{n\rightarrow\infty}|a_n|^{-1/n}\geq \rho\r\}\ ,
$$
of analytic functions (with coefficients in $K$) over the open disk
$$
D=D(\xi,\rho^-)=\{x\in K :|\xi-x|<\rho\}
$$
of center $\xi\in{\mathbb A}^1_K$ and radius $\rho\in\R$, $\rho>0$,
and the field of meromorphic functions ${\mathcal M}_D=Frac({\mathcal A}_D)$ over $D$.
Sometimes we will write ${\mathcal A}_{D,K}$ (resp. ${\mathcal M}_{D,K}$)
to stress the fact that we are considering analytic (resp. meromorphic) functions
with coefficients in $K$.
\par
If $D$ is $q$-invariant
(i.e. if $D$ is invariant for the isometry $x\longmapsto qx$, or,
equivalently, if $|(q-1)\xi|<\rho$), it makes sense to consider
a \emph{$q$-difference systems (of order $\mu$, with meromorphic coefficients over $D$, defined over $K$)}:
$$
Y(qx)=A(x)Y(x)\ ,\
\hbox{with }A(x)\in Gl_\mu({\mathcal M}_{D,K}).
\leqno{(\mathcal S)_q}$$
The main purpose of this chapter is to introduce the basic properties of $q$-difference systems.
First we study the properties of the
$q$-difference algebra of analytic functions and of $q$-difference operators
with respect to Gauss norms.
Then we construct analytic solutions of $q$-difference systems, when they exist.
Finally, we define the notions of generic point and generic radius of convergence.

\subsection{The $q$-difference algebra of analytic functions over an open disk}
\label{subsec:anfunction}

Let the disk $D=D(\xi,\rho^-)$ be $q$-invariant.
Then the \emph{$q$-difference operator}
$$
\sgq:f(x)\longmapsto f(qx)
$$
is a $K$-algebra isomorphism of ${\mathcal A}_D$:
we say that ${\mathcal A}_D$ is a \emph{$q$-difference algebra}.

\medskip
One can also define a \emph{$q$-derivation}
$$
\dq(f)(x)=\frac{f(qx)-f(x)}{(q-1)x}\ ,
$$
satisfying the twisted Leibniz rule
$$
\dq(fg)=\sgq(f)\dq(g)+\dq(f)g\ .
$$

\begin{lemma}
The operator $\dq$ acts over ${\mathcal A}_D$.
\end{lemma}

\begin{proof}
Let $0\not\in D$. Then by definition
$\dq(f)\in{\mathcal A}_D$ for all $f\in{\mathcal A}_D$.
\par
If $0\in D$, we can suppose $\xi=0$.
Observe that
\begin{equation}\label{eq:coeff}
\dq x^n=\l(1+q+\dots+q^{n-1}\r) x^{n-1}
\ \forall n\geq 1.
\end{equation}
Then for any $\sum_{n\geq 0}a_nx^n\in{\mathcal A}_D$ we have:
$$
\dq\l(\sum_{n\geq 0}a_nx^n\r)=
\sum_{n\geq 1}\l(1+q+\dots+q^{n-1}\r) a_nx^{n-1}\in{\mathcal A}_D\ .
$$
\end{proof}

\medskip
Motivated by (\ref{eq:coeff}), we recall the classical definition of
\emph{$q$-factorials} and
\emph{$q$-binomial coefficients}, namely for any pair of integers $n\geq i\geq 0$ we set:
$$
 \begin{array}{cl}
&[0]_q=0\;,\;\nq=1+q+\dots+q^{n-1}=\frac{1-q^n}{1-q}\;,\\ \\
&[0]_q^!=1\;,[n]_q^!=1_q\cdots\nq\;,\\ \\
&\ds{n\choose 0}_q={n\choose n}_q=1\;,\;
       \ds{n\choose i}_q=\frac{[n]_q^!}{[n-i]_q^![i]_q^!}=\frac{\nq [n-1]_q\cdots[n-i+1]_q}{[i]_q^!}\;.\\
 \end{array}
$$
They satisfy the relation
$$
{n\choose i}_q
={n-1\choose i-1}_q+{n-1\choose i}_q q^i
={n-1\choose i-1}_q q^{n-i}+{n-1\choose i}_q\ ,
\ \hbox{for $n\geq i\geq 1$,}
$$
and
\begin{equation}\label{eq:binomio}
(1-x)(1-qx)\cdots(1-q^{n-1}x)=\sum_{i=0}^n(-1)^i{n\choose i}_qq^{\frac{i(i-1)}{2}}x^i\ .
\end{equation}
One verifies directly the following basic properties of $\sgq$ and $\dq$:

\begin{lemma}\label{lemma:leibniz}
For any pair integers $n,i\geq 1$ and any $f,g\in{\mathcal A}_D$ we have:

\begin{eqnarray}
&\hskip 20 pt\hbox to 11 true cm{$\ds\frac{\dq^n}{[n]_q^!}x^i=
    \begin{cases}
    {n\choose i}_qx^{n-i},&\hbox{ if }n\geq i,\\
    0,&\hbox{ otherwise;}
    \end{cases}$\hfill}\label{eq:leibnizuno}\\
&\hskip 20 pt\hbox to 11 true cm{$\ds\dq^n(fg)(x)=
    \sum_{j=0}^n {n\choose j}_q\dq^{n-j}(f)(q^jx)\dq^j(g)(x)$;\hfill}\label{eq:leibnizdue}\\
&\hskip 20 pt\hbox to 11 true cm{$\ds\sgq^n=\sum_{j=0}^n{n\choose j}_q(q-1)^j q^{\frac{j(j-1)}{2}}x^j\dq^j
    =\sum_{j=0}^n\prod_{i=0}^{j-1}(q^n-q^i)x^j
    \frac{\dq^j}{[j]_q^!}$;\hfill}\label{eq:leibniztre}\\
&\hskip 20 pt\hbox to 11 true cm{$\ds\dq^n
    =\frac{\prod_{i=0}^{j-1}(\sgq-q^i)}{(q-1)^nq^{\frac{n(n-1)}{2}}x^n}
    =\frac{(-1)^n}{(q-1)^nx^n}
    \sum_{j=0}^n (-1)^j{n\choose j}_{q^{-1}}q^{-\frac{j(j-1)}{2}}\sgq^j$;\hfill}
    \label{eq:leibnizquattro}\\
&\hskip 20 pt\hbox to 11 true cm{$\ds(x\dq)^n =\frac{1}{(q-1)^n}
    \sum_{j=0}^n (-1)^j{n\choose j}\sgq^j$.\hfill}
    \label{eq:leibnizcinque}
\end{eqnarray}
\end{lemma}

\begin{proof}
The first four formulas are proved in \cite[(1.1.8) through (1.1.10)]{DV}.
The proof of (\ref{eq:leibnizcinque}) is straightforward:
$$
(x\dq)^n = \l(\frac{\sgq-1}{q-1}\r)^n
= \frac{1}{(q-1)^n}\sum_{i=0}^n(-1)^i{n\choose i}\sgq^i\ .
$$
\end{proof}

\begin{parag}\label{parag:mult}
The topological basis $\l((x-\xi)^n\r)_{n\geq 0}$ of ${\mathcal A}_D$ is not adapted to
study the action of the $q$-derivation over ${\mathcal A}_D$ as the relation
$$
\dq (x-\xi)^n=\frac{(qx-\xi)^n-(x-\xi)^n}{(q-1)x}
$$
clearly shows. So, rather than $(x-\xi)^n$, one classically consider the polynomials
$$
\begin{array}{l}
(x-\xi)_{q,0}=1\;,\\
(x-\xi)_{q,n}=(x-\xi)(x-q\xi)\cdots(x-q^{n-1}\xi)\hbox{, for any integer $n\geq 1$,}
\end{array}
$$
which satisfy the formula (cf. (\ref{eq:coeff}))
$$
\dq (x-\xi)_{q,n}=[n]_q(x-\xi)_{q,n-1}\ .
$$
Therefore we are naturally led to introduce the $q$-difference
algebra
$$
K\{x-\xi\}_{q,\rho}=
\l\{\sum_{n\geq 0}a_n(x-\xi)_{q,n}:a_n\in K,
\liminf_{n\rightarrow\infty}|a_n|^{-1/n}\geq\rho\r\}\ .
$$
If $\xi=0$, the $K$-algebra $K\{x-\xi\}_{q,\rho}$ obviously conincides with
the ring of analytic functions over the open disk of center $0$ and
radius $\rho$ and its structure is well-known.
If $\xi\neq 0$, there are different equivalent ways to define the multiplication of $K\{x-\xi\}_{q,\rho}$.
Let $f(x)=\sum_{n\geq 0}f_n(x-\xi)_n,g(x)=\sum_{n\geq 0}g_n(x-\xi)_n
\in K\{x-\xi\}_{q,\rho}$. Then for any nonnegative integers $i,n$ we
have
\begin{equation*}
\frac{\dq^n}{[n]_q^!}(f)(\xi)=f_n
\hbox{ and }
f(q^i\xi)=\sum_{n=0}^i(q^i-1)(q^1-q)\cdots(q^i-q^{n-1})\xi^n f_n\ .
\end{equation*}
Moreover $f(x)$ is uniquely determined by the sequence
$f(q^i\xi)$, $i\geq 0$.
Analogous relations hold for $g(x)$. Hence the product $f(x)g(x)$
is the element of $K\{x-\xi\}_{q,\rho}$ associated to the
sequence $f(q^i\xi)g(q^ix)$, namely, if
$f(x)g(x)=\sum_{n\geq 0}h_n(x-\xi)_n$, we deduce from
(\ref{eq:leibnizquattro}) that
$$
h_n=\frac{\dq^n}{[n]_q^!}(fg)(\xi)
=\frac{(-1)^n}{(q-1)^n[n]_q^!\xi^n}
    \sum_{j=0}^n (-1)^j{n\choose j}_{q^{-1}}
    q^{-\frac{j(j-1)}{2}}f(q^j\xi)g(q^j\xi)\ .
$$
The twisted Leibniz Formula (\ref{eq:leibnizdue}) gives another natural way of defining
the coefficient $h_n$, in fact
$$
h_n=\frac{\dq^n}{[n]_q^!}(fg)(\xi)=
    \sum_{j=0}^n
    \sum_{h=j}^n {h\choose j}_q f_h(q^{n-j}-1)_{h-j}\xi^{h-j}g_{n-j}\ .
$$
In particular for any pair of positive integers $l,k$ the formula above
specializes to
$$
(x-\xi)_l(x-\xi)_k=\sum_{n=0}^{l+k}
{l\choose n-k}_q(q^k-1)_{l+k-n}\xi^{l+k-n}(x-\xi)_n\ .
$$
\end{parag}

The following proposition states that
the natural map
\begin{equation}\label{eq:defT}
T_{q,\xi}:
f(x)\longmapsto\ds\sum_{n\geq 0}
\frac{d_q^nf}{[n]_q^!}(\xi)(x-\xi)_{q,n}
\end{equation}
defines an isomorphism of
$q$-difference algebras
(i.e. an isomorphism of $K$-algebras commuting to the action of $\dq$)
between ${\mathcal A}_D$ and $K\{x-\xi\}_{q,\rho}$.
We will call the map $T_{q,\xi}$ \emph{$q$-expansion}
or \emph{twisted Taylor formula}.

\begin{prop}\label{prop:taylor}
Let $D=D(\xi,\rho^-)$ be a $q$-invariant open disk.
The map
$$
\begin{array}{rccl}
T_{q,\xi}:&{\mathcal A}_D&\longrightarrow & K\{x-\xi\}_{q,\rho}\\ \\
          &f(x)          &\longmapsto\ds  &\sum_{n\geq 0}
               \ds \frac{d_q^nf}{[n]_q^!}(\xi)(x-\xi)_{q,n}
\end{array}
$$
is a $q$-difference algebras isomorphism.
Moreover, for all $f\in{\mathcal A}_D$, the series $T_{q,\xi}(f)(x)$
converges uniformly to $f(x)$ over any closed disk $D(\xi,\eta^+)$, with
$0<\eta<\rho$.
\end{prop}

The proof of (\ref{prop:taylor}) can be found in the Appendix
(cf. \S\ref{subsec:qtaylorconnesso}).
As a corollary we obtain the more useful statement:

\begin{cor}
Let $f(x)=\sum_{n\geq 0}a_n(x-\xi)_{q,n}$ be a series such that $a_n\in K$
and let $\ds\rho=\liminf_{n\rightarrow\infty}|a_n|^{-1/n}$.
Then $f(x)$ converges uniformly over any closed disk $D(\xi,\eta^+)$, with
$0<\eta<\rho$, to an analytic function
if and only if $\rho>|(q-1)\xi|$.
\end{cor}

\begin{proof}
If $\rho>|(q-1)\xi|$, the series $f(x)$ converges by (\ref{prop:taylor}).
\par
Suppose that $\rho<|(q-1)\xi|$. Let
$n_0$ be the smallest positive integer such that
$|(q^{n_0}-1)\xi|\leq \rho$ and let $\veps$ be a real positive
number such that
$$
\sup_{i=0,\dots,n_0-1}|(q^i-1)\xi|>\rho+\veps\ .
$$
Then for any $k\in\Z_{>0}$ we have
$$
|(q^{kn_0}-1)\xi|\leq \sup_{i=0,\dots,k-1}|(q^{n_0}-1)q^{in_0}\xi|\leq \rho
$$
and for any $r=1,\dots,n_0-1$ we have
$$
|(q^{kn_0+r}-1)\xi|=\sup\l(|(q^{kn_0}-1)\xi|,|(q^r-1)q^{kn_0}\xi|\r)>\rho+\veps\ .
$$
Therefore for $x_0\in D(\xi,\rho^-)$ we obtain
$$
\begin{array}{rl}
|(x_0-\xi)_{n,q}|
&=\ds\prod_{i=0}^{n-1}|(x_0-\xi)+\xi(1-q^i)|\\
&>(\rho+\veps)^{n-\l[\frac{n-1}{n_0}\r]-1}\ds\prod_{0\leq i\leq\l[\frac{n-1}{n_0}\r]}|x_0-\xi+\xi(1-q^{in_0})|\ .
\end{array}
$$
We conclude that
$$
\sup_{|x_0-\xi|<\rho}|(x_0-\xi)_{n,q}|
\geq (\rho+\veps)^{n-\l[\frac{n-1}{n_0}\r]-1}\rho^{\l[\frac{n-1}{n_0}\r]+1}
>\rho^n
$$
and hence that $\sum_{n\geq 0}a_n(x_0-\xi)_{q,n}$ does not converge over $D(\xi,\rho^-)$.
\par
Let $\rho=|(q-1)\xi|$. Suppose that $\sum_{n\geq 0}a_n(x-\xi)_{q,n}$ is convergent over
$D(\xi,\rho^-)$. Then the series
$$
f(qx)=a_0+(qx-\xi)q\sum_{n\geq 1}a_n(x-\xi)_{q,n-1}
$$
must also converge. By induction we conclude that
for any $x_0\in D(\xi,\rho^-)$ and any integer $n\geq 0$, the sum
$f(q^nx_0)$ is convergent. Hence $f(x)$ converges over a bigger disk than the closed disk of radius $\rho$,
which means that $\ds\liminf_{n\rightarrow\infty}|a_n|^{-1/n}>\rho$.
\end{proof}

\subsection{Gauss norms and $q$-difference operators}
\label{subsec:gaussnorm}

Let $D=D(\xi,\rho^-)$.
We recall that ${\mathcal A}_D$ comes equipped with
a family of non archimedean norms $|~|_\xi(R)$, the so-called Gauss norms,
(cf. for instance \cite[\S6, 1.4]{ROB})
$$
\l\|\sum_{n\geq 0}f_n(x-\xi)^n\r\|_\xi(R)
=\sup_{n\geq 0}|f_n|R^n\ ,
$$
defined for any $R\in (0,\rho)$ and any $\sum_{n\geq 0}f_n(x-\xi)^n\in{\mathcal A}_D$.
It follows by Gauss lemma that they are multiplicative norms.
If moreover $R\in| K |$, then (cf. \cite[IV, 1.1]{DGS})
$$
\l\|\sum_{n\geq 0}f_n(x-\xi)^n\r\|_\xi(R)
=\sup_{x_0\in K ,~|x_0-\xi|< R}\l|\sum_{n\geq 0}f_n(x_0-\xi)^n\r|\ .
$$
The norm $\|~\|_\xi(R)$ plays a central role in the study of $p$-adic
differential equation as well as of $p$-adic $q$-difference equations; therefore
it is crucial to calculate the norm of the $K$-linear operator
$\frac{d_q^k}{[k]_q^!}$ with respect to $\|~\|_\xi(R)$
as well as to be able to calculate the norm of $f\in {\mathcal A}_D$
by looking at its $q$-expansion:

\begin{prop}
Let $f(x)\in{\mathcal A}_{D(\xi,\rho^-)}$. If $|(1-q)\xi|\leq R<\rho$,
the norm $\|~\|_\xi(R)$ satisfies
\begin{equation}\label{eq:normaqder}
    \l\|\frac{d_q^k}{[k]_q^!}f(x)\r\|_\xi(R)
    \leq\frac{1}{R^k}\l\|f(x)\r\|_\xi(R)\
\hbox{ for all }k\in\Z_{\geq 0};
\end{equation}
and
\begin{equation}\label{eq:supnorma}
\|f(x)\|_\xi(R)=\sup_{n\geq 0}\l|\frac{d_q^n(f)}{[n]_q^!}(\xi)\r|R^n\ .
\end{equation}
\end{prop}

\medskip\noindent
{\bf Remark.}
Inequality (\ref{eq:normaqder}) generalizes
\cite[4.2.1]{DV}, where we considered the case $\xi=0$ and $R=1$.

\begin{proof}
Let us prove (\ref{eq:normaqder}).
Since
$$
\frac{d_q^k}{[k]_q^!}f(x)=\sum_{n\geq k}f_n
       \frac{d_q^k}{[k]_q^!}(x-\xi)^n\ ,
$$
it is enough to prove that
$$
\l\|\frac{d_q^k}{[k]_q^!}(x-\xi)^n\r\|_\xi(R)
\leq \frac{1}{R^k}\l\|(x-\xi)^n\r\|_\xi(R)=R^{n-k}\ ,
\hbox{ for any $n\geq k$.}
$$
We proceed by double induction over $k,n$.
\par
Let $k=1$.
If $n$ is an odd positive integer the inequality immediately
follows from
$$
\dq(x-\xi)^n=\frac{(qx-\xi)^n-(x-\xi)^n}{(q-1)x}
=\sum_{i=0}^{n-1}(qx-\xi)^{n-1-i}(x-\xi)^i\ .
$$
If $n$ is an even positive integer we reduce to the previous case
by observing that
$$
\begin{array}{rcl}
\dq(x-\xi)^n
&=&\ds\frac{(qx-\xi)^n-(x-\xi)^n}{(q-1)x}\\
&=&\ds\frac{(qx-\xi)^{n/2}-(x-\xi)^{n/2}}{(q-1)x}
       \l((qx-\xi)^{n/2}+(x-\xi)^{n/2}\r)\\
&=&\ds\dq(x-\xi)^{n/2}\l((qx-\xi)^{n/2}+(x-\xi)^{n/2}\r)\ .
\end{array}
$$
Now let $k>1$.
The inequality is clear for $n=k$.
It follows from the twisted Leibniz Formula
(\ref{eq:leibnizdue}) that for any $n\geq k$ we have
$$
\begin{array}{rcl}
\lefteqn{\l\|\frac{d_q^k}{[k]_q^!}(x-\xi)^{n+1}\r\|_\xi(R)}\\
&=&\l\|(q^kx-\xi)\frac{d_q^k}{[k]_q^!}(x-\xi)^n+
      \frac{d_q^{k-1}}{[k-1]_q^!}(x-\xi)^n\r\|_\xi(R)\\
&\leq&\sup\l(\l\|(q^kx-\xi)\frac{d_q^k}{[k]_q^!}
      (x-\xi)^n\r\|_\xi(R),
      \l\|\frac{d_q^{k-1}}{[k-1]_q^!}
      (x-\xi)^n\r\|_\xi(R)\r)\\
&\leq& R^{n+1-k}\ ,
\end{array}
$$
which completes the proof of (\ref{eq:normaqder}).
\par
Clearly (\ref{eq:normaqder}) implies that
$$
\|f(x)\|_\xi(R)\geq\sup_{n\geq 0}\l|\frac{d_q^n(f)}{[n]_q^!}(\xi)\r|R^n\ .
$$
So it is enough to prove the opposite inequality to obtain (\ref{eq:supnorma}).
By (\ref{prop:taylor}) we have $f(x_0)=T_{q,\xi}(f)(x_0)$, for any $x_0\in D$.
Then it is enough to remark that
$$
\l\|(x-\xi)_{q,n}\r\|_\xi(R)=
\l\|\prod_{i=0}^{n-1}\Big[(x-\xi)+\xi(1-q^i)\Big]\r\|_\xi(R)
\leq R^n\ \hbox{ for any integer $n\geq 0$,}
$$
to conclude that
$$
\|f(x)\|_\xi(R)
\leq\sup_{n\geq 0}\l|\frac{d_q^n(f)}{[n]_q^!}(\xi)\r|\l\|(x-\xi)_{q,n}\r\|_\xi(R)
\leq\sup_{n\geq 0}\l|\frac{d_q^n(f)}{[n]_q^!}(\xi)\r|R^n\ .
$$
This finishes the proof.
\end{proof}

The following lemma will be useful in (\ref{parag:deformation}), where we
will consider some properties of families of $q$-difference systems deforming a differential system.

\begin{lemma}\label{lemma:stimaqdef}
Let $f(x)=\sum_{n\geq 0}a_n(x-\xi)^n$ and $g(x)=\sum_{n\geq 0}b_n(x-\xi)_{q,n}$
be two analytic bounded functions over $D(\xi,1^-)$, with $|\xi|\leq 1$,
and let $0<|1-q|\leq\veps$.
\par
Suppose $\l\|f-g\r\|_\xi(1)\leq \veps$.
Then  $\l\|\frac{df}{dx}-\dq(g)\r\|_\xi(1)\leq \veps$
\end{lemma}

\begin{proof}
Notice that for any positive integer $n$ we have
$$
|[n]_q-n|=|(q-1)+\dots+(q^{n-1}-1)|\leq \sup_{i=1,\dots,n-1}|q^i-1|\leq |q-1|<\veps
$$
and
$$
\l\|(x-\xi)^n\r\|_\xi(1)=\l\|(x-\xi)_{q,n}\r\|_\xi(1)=1\ .
$$
We can assume that both $\|f\|_\xi(1)\leq 1$ and $\|g\|_\xi(1)\leq 1$.
Therefore one conclude that
$$
\begin{array}{rcl}
\l\|\frac{df}{dx}-\dq(g)\r\|_\xi(1)
&=&\l\|\sum_{n\geq 1}a_n n!(x-\xi)^{n-1}-b_n[n]_q^!(x-\xi)_{q,n-1}\r\|_\xi(1)\\ \\
&\leq&\sup_{n\geq 1}\Big(\l|a_n n!-b_n n!\r|,
    \l|b_n n!-b_n [n]_q^!\r|,\\ \\
&&\l|b_n [n]_q^!\r|\l\|(x-\xi)^{n-1}-(x-\xi)_{q,n-1}\r\|_\xi(1)\Big)\\ \\
&\leq&\sup\l(|n!|\veps, |(1-q)b_n|,\l|b_n [n]_q^!(1-q)\r| \r)\\ \\
&\leq&\sup(\veps,|1-q|)\leq\veps\ .
\end{array}
$$
\end{proof}

\subsection{Analytic solutions of $q$-difference systems}
\label{subsec:analyticsolution}

Let us consider a $q$-difference system of order $\mu$
$$
Y(qx)=A(x)Y(x)\ ,
\leqno{(\mathcal S)_q}$$
whose coefficients are meromorphic functions over a $q$-invariant open disk $D(\xi,\rho^-)$.
The system $(\mathcal S)_q$ can be rewritten in the form
$$
  \dq(Y)(x)=G(x)Y(x), \hbox{ with } G(x)=\frac{A(x)-{\mathbb I}_\mu}{(q-1)x}\ ,
  \leqno{({\mathcal S})_q^\p}
$$
where ${\mathbb I}_\mu$ is the identity matrix of order $\mu$.
One can iterate $({\mathcal S})_q^\p$ obtaining
$$
\dq^n(Y)(x)=G_n(x)Y(x)\ ,
$$
with
\begin{equation}\label{eq:ricorsione}
\begin{array}{l}
\hbox{$G_0(x)={\mathbb I}_\mu$, $G_1(x)=G(x)$ and } \\
G_{n+1}(x)=G_n(qx)G(x)+\dq(G_n)(x), \hbox{ for any integer $n\geq 1$.}
\end{array}
\end{equation}
If $G(x)$ does not have any pole in $q^\N\xi=\{\xi,q\xi,q^2\xi,\dots\}$,
it can be identified, by considering its $q$-expansion, with a matrix with entries in the ring
$$
K\[[x-\xi\]]_q=\l\{\sum_{n\geq 0}a_n(x-\xi)_{q,n}:a_n\in K\r\}\ ,
$$
whose multiplication is defined as in (\ref{parag:mult}).
Then a formal solution matrix of $(\mathcal S)_q$, or equivalently of $(\mathcal S)_q^\p$,
at $\xi$ is given by
\begin{equation}\label{eq:analyticsolution}
Y(\xi,x)=\sum_{n\geq 0}\frac{G_n(\xi)}{[n]_q^!}(x-\xi)_{q,n}\ .
\end{equation}
The fact that $Y(\xi,\xi)={\mathbb I}_\mu$ does not allow us to conclude that
$Y(\xi,x)$ is an invertible matrix, since $K\[[x-\xi\]]_q$ is not a local ring
(cf. Appendix \ref{sec:taylor}, \S\ref{subsec:taylor}).
Actually we need a stronger assumption:

\begin{lemma}\label{lemma:solformali}
The system $(\mathcal S)_q$ has a formal solution matrix in $Gl_\mu(K\[[x-\xi\]]_q)$
if and only if the matrix $A(x)$ does not have any poles in $q^\N\xi$ and
$\det A(x)$ does not have any zeros in $q^\N\xi$.
\end{lemma}

\begin{rmk}\hfill
\par\smallskip\noindent
1) If the conditions of the lemma above are verified,
then (\ref{eq:analyticsolution}) is the only solution of $(\mathcal S)_q$
in $Gl_\mu(K\[[x-\xi\]]_q)$ such that $Y(\xi,\xi)={\mathbb I}_\mu$ and all
other solution matrices of  $(\mathcal S)_q$
in $Gl_\mu(K\[[x-\xi\]]_q)$ are obtained by multiplying $Y(\xi,x)$ on the right by an element of
$Gl_\mu(K)$.
\par\smallskip\noindent
2) Observe that, if $({\mathcal S})_q$ has a solution matrix $Y(x)\in Gl_\mu({\mathcal A}_D)$
over a $q$-invariant disk $D$,
the matrix $A(x)=Y(qx)Y(x)^{-1}$ is an element of $Gl_\mu({\mathcal A}_D)$. Hence neither $A(x)$
has a pole in $q^\Z \xi$ nor $\det A(x)$ has a zero in $q^\Z\xi$.
It follows by the statement above that $({\mathcal S})_q$ can have a formal solution in
$Gl_\mu(K\[[x-\xi\]]_q)$ which is not the $q$-expansion of an analytic solution.
\par\smallskip\noindent
3) Suppose that $Y(\xi,x)\in Gl_\mu(K\[[x-\xi\]]_q)$
is the $q$-expansion of an analytic solution of $({\mathcal S})_q$
converging over $D(\xi,\rho^-)$ and let $|\zeta-\xi|<\rho$. Then necessarily we have
\begin{equation}\label{eq:solformali}
Y(\zeta,x)=Y(\xi,x)Y(\xi,\zeta)^{-1}\in Gl_\mu(K\[[x-\zeta\]]_q)\ ,
\end{equation}
since both matrices are analytic solution matrix of $({\mathcal S})_q$ at $\zeta$, of
maximal rank, having value ${\mathbb I}_\mu$ at $\zeta$.
\end{rmk}

\begin{proof}
By the remark above, the system $(\mathcal S)_q$ has a formal solution matrix in the ring $Gl_\mu(K\[[x-\xi\]]_q)$
if and only if $Y(\xi,x)$ is in $Gl_\mu(K\[[x-\xi\]]_q)$.
\par
For any non negative integer $k$ it makes sense to evaluate $Y(\xi,x)$ at $q^k\xi$:
$$
Y(\xi,q^k\xi)=\sum_{n\geq 0}\frac{G_n(\xi)}{[n]_q^!}\xi^n(q^k-1)_{q,n}\ ,
$$
since the sum on the right hand side is actually finite.
Of course there are precise relations between the sequences
$\l(\frac{G_n(\xi)}{[n]_q^!}\r)_{k\geq 0}$
and $(Y(\xi,q^k\xi))_{k\geq 0}$ that one can easily deduce by
(\ref{eq:leibniztre}) and (\ref{eq:leibnizquattro}).
It turns out that an element $Y(x)\in M_{\mu\times\mu}(K\[[x-\xi\]]_q)$
is uniquely determined by $(Y(q^k\xi))_{k\geq 0}$. Therefore
$Y(\xi,x)\in Gl_\mu(K\[[x-\xi\]]_q)$ if and only if
$Y(\xi,q^k\xi)\in Gl_\mu(K)$ for all $k\geq 0$, the inverse of $Y(\xi,x)$ being
the element of $Gl_\mu(K\[[x-\xi\]]_q)$ associated to the data $\l(Y(\xi,q^k\xi)^{-1}\r)_{k\geq 0}$.
\par
To conclude it is enough  to observe that
$$
\begin{array}{rcl}
Y(\xi,q^k\xi)&=&A(q^{k-1}\xi)A(q^{k-2}\xi)\cdots A(\xi)Y(\xi,\xi)\\
&=&A(q^{k-1}\xi)A(q^{k-2}\xi)\cdots A(\xi)\ .
\end{array}
$$
\end{proof}

In the next corollaries we give some sufficient conditions for having
a fundamental analytic solution matrix, i.e.
an invertible solution matrix $Y(x)$ such that $Y(x)$ and $Y(x)^{-1}$
have analytic coefficients over a convenient $q$-invariant open disk.
It is just a partial result and
we will reconsider the problem of the existence of analytic solutions
for $({\mathcal S})_q$ in the next sections.

\medskip
From (\ref{prop:taylor}) one immediately obtains:

\begin{cor}\label{cor:solanalytic}
The system $(\mathcal S)_q$ has a fundamental analytic solution matrix at $\xi$
if and only if
\begin{trivlist}
\item - the matrix $A(x)$ does not have any poles in $q^\N\xi$,
\item - $\det A(x)$ does not have any zeros in $q^\N\xi$,
\item - $\ds\limsup_{n\rightarrow\infty} \l|\frac{G_n(\xi)}{[n]_q^!}\r|^{1/n}<|(q-1)\xi|^{-1}$.
\end{trivlist}
\end{cor}

\medskip
Before stating the following result we need
to introduce the number $\pi_q$, which plays a role analogous to the $\pi$ of Dwork for
$p$-adic differential equations. We recall that $\pi$ is an element
of $ K $ such that $\pi^{p-1}=-p$.

\begin{notation}\label{defn:pi}
We fix an element $\pi_q$ of $K$ such that
$\ds\lim_{n\rightarrow\infty} |[n]_q^!|^{1/n}=|\pi_q|$.
\end{notation}

\begin{rmk}
The condition above is somehow a minimal require for $\pi_q$, and
actually it does not define it uniquely. Further
developments of the theory point out extra conditions that
will determine more precisely the choice of $\pi_q$ (cf. \cite{Villetaneuse}).
\par
A possible choice for $\pi_q$ is the following.
Let $\kappa$ be the smallest positive integer such that $|1-q^\kappa|<|\pi|$. Then
$\pi_q$ can be chosen to be an element of $ K $ defined by
$$
\pi_q^{p-1}=-[p]_q,\hbox{ if } \kappa=1, \hbox{ and }
\pi_q^\kappa=[\kappa]_q\pi_{q^\kappa},\hbox{ if } \kappa>1.
$$
Then it follows from \cite[4.1.1]{DV} that
$\pi_q$ has the desired property.
In particular $|\pi_{q^\kappa}|=|\pi|$.
\end{rmk}

\begin{cor}
Let $\rho\leq 1$, $\rho|\pi_q|>|(1-q)\xi|$ and $D=D(\xi,\rho^-)$.
Suppose we are given a square matrix $G(x)$ analytic over $D$ such that
$$
\sup_{x\in D}\l|G(x)\r|\leq\frac{1}{\rho}
$$
and that the determinant of $A(x)=(q-1)xG(x)+{\mathbb I}_\mu$
does not have any zeros in $q^\N\xi$.
Then $Y(qx)=A(x)Y(x)$ has an analytic fundamental solution at $\xi$.
\end{cor}

\begin{proof}
It follows by (\ref{eq:normaqder}) and (\ref{eq:ricorsione}) that
$$
\l|G_n(\xi)\r|\leq\frac{1}{\rho}\sup\l(1,\l\|G_{n-1}(\xi)\r\|\r)\leq \frac{1}{\rho^{n}}\ ,
\hbox{ for any $n\geq 1$,}
$$
which implies that
$$
\limsup_{n\rightarrow\infty} \l|\frac{G_n(\xi)}{[n]_q^!}\r|^{1/n}=
\frac{1}{|\pi_q|}\limsup_{n\rightarrow\infty} \l|G_n(\xi)\r|^{1/n}\leq \frac{1}{\rho|\pi_q|}\ .
$$
Since $\rho|\pi_q|>|(1-q)\xi|$, the matrix $Y(x,\xi)$ in (\ref{eq:analyticsolution})
is the $q$-expansion of an analytic fundamental solution.
\end{proof}

\begin{parag}
{\bf Iteration of $({\mathcal S})_q$ and existence of analytic solutions.}

\par\medskip\noindent
In the rest of the paper we will often assume that $(\mathcal S)_q$ has
an analytic fundamental solution at some point $\zeta$ or that
$|1-q|$ is smaller than some constant: this is not always true.
\par
Sometimes one can
easily reduce to the case of having a fundamental analytic solution
by iterating the $q$-difference system.
In the same way one can reduce to the case of a $q\in K$ such that $|1-q|<<1$.
\par
Let us analyze the situation in detail.
Suppose that $(\mathcal S)_q$ does not have
an analytic fundamental solution at $\zeta$.
Then it may happen that there exists $n_0>1$ such that the system
$$
Y(q^{n_0}x)=A_{n_0}(x)Y(x)\ ,\
\hbox{with }A_{n_0}(x)=A(q^{n_0-1}x)A(q^{n_0-2}x)\cdots A(x),
\leqno{(\mathcal S)_{q^{n_0}}}
$$
obtained from $(\mathcal S)_q$ by iteration, has a fundamental analytic solution $Y(x)$
over a $q^{n_0}$-invariant open disk $D(\zeta, \eta^-)$.
If $n_0$ is the smallest positive integer having this property,
then one can construct a fundamental solution $F(x)$ of
$(\mathcal S)_q$, analytic over the non-connected $q$-invariant analytic domain
\begin{equation}\label{eq:primoqdisco}
D(\zeta, \eta^-)\cup D(q\zeta, \eta^-)\cup\dots\cup D(q^{n_0-1}\zeta, \eta^-)\
\end{equation}
by setting
$$
\begin{array}{l}
F(q^ix)=A_i(x)Y(x)=A(q^{i-1}x)A(q^{i-2}x)A(q^{i-1}x)\cdots A(x)Y(x)\ ,\\
\hbox{for any $x\in D(\zeta, \eta^-)$ and any $i=0,\dots,n_0-1$.}
\end{array}
$$
Observe that the restriction of $F(x)$ to $D(q^i\zeta, \eta^-)$, for any $i\in\Z$, is an analytic
fundamental solution of $(\mathcal S)_{q^{n_0}}$.
\par
So, if the $q$-difference system has an analytic solution over a non-connected analytic domain as above,
it is enough to consider a system obtained by iteration to reduce to the case of a system having
an analytic solution over a $q$-invariant open disk.
In the appendix we will consider the $q$-expansion of analytic functions over
non connected domain of the form
(\ref{eq:primoqdisco}).
\par
The same trick allows us to reduce to the case of a small $|1-q|$,
knowing that $\ds\liminf_{n\rightarrow\infty}|1-q^n|=0$.
\end{parag}

\begin{parag}
{\bf Removing apparent and trivial singularities.}

\par\medskip\noindent
In this subsection we will consider $q$-difference system having
meromorphic solutions or analytic solutions with meromorphic inverse.
Our purpose is to explain how to reduce
by gauge transformation to the assumption of having
an analytic fundamental solution.

\medskip
Once again we consider a $q$-difference system
$$
Y(qx)=A(x)Y(x)\
\leqno{(\mathcal S)_q}
$$
with meromorphic coefficients over a $q$-invariant disk $D=D(\xi,\rho^-)$, defined over $K$.
For any matrix $F(x)\in Gl_\mu({\mathcal M}_D)$, the matrix
$Z(x)=F(x)Y(x)$ is a solution to
\begin{equation}\label{eq:gauge}
Z(qx)=A_{[F]}(x)Z(x)\ ,
\hbox{ with }A_{[F]}(x)=F(qx)A(x)F(x)^{-1}.
\end{equation}
The matrix $F(x)$ is usually called a \emph{meromorphic gauge transformation matrix}.
Notice that
\begin{equation}\label{eq:dqgauge}
G_{[F]}(x)
=\frac{A_{[F]}(x)-{\mathbb I}_\mu}{(q-1)x}=F(qx)G(x)F(x)^{-1}+\dq(F)(x)F(x)^{-1}\ .
\end{equation}

\medskip
Following the classical terminology of $p$-adic differential equations
(cf. for instance \cite[page 172]{DGS})
we give the definition:

\begin{defn}
We say that $q^\Z\xi\subset D$ is an
\emph{ordinary $q$-orbit ({\rm resp.} trivial singularity, apparent singularity) for $({\mathcal S})_q$},
if $({\mathcal S})_q$ has a solution in $Gl_\mu({\mathcal A}_{D^\p})$
(resp. $Gl_\mu({\mathcal M}_{D^\p})\cap M_{\mu\times\mu}({\mathcal A}_{D^\p})$,
$Gl_\mu({\mathcal M}_{D^\p})$),
where $D^\p\subset D$ is a $q$-invariant analytic domain of the form
(\ref{eq:primoqdisco}) containing $\xi$.
\end{defn}

\begin{rmk}
In the sequel we will informally say that the system $({\mathcal S})_q$ has ``\emph{at worst a finite number of
apparent singularities in $D$}'' to mean that $({\mathcal S})_q$ has only ordinary $q$-orbits in $D$,
apart from a finite number of apparent singularities.
\end{rmk}

\medskip
The following statement is a $q$-analogue of the Frobenius-Christol device
(cf. \cite[II, \S8]{Cfrobenius})
to remove apparent and trivial singularities over
$$
D^\times=\begin{cases}
D& \hbox{if $0\not\in D$,}\\
D\smallsetminus\{0\} & \hbox{otherwise.}
\end{cases}
$$

\begin{prop}\label{prop:rem}
We assume that
\begin{equation}\label{eq:rem}
\hbox{the system $({\mathcal S})_q$ has at worst a finite number of apparent singularities in $D^\times$.}
\end{equation}
Then there exists $H(x)\in Gl_\mu(K(x))$ such that the $q$-difference system
$Y(qx)=A_{[H]}(x)Y(x)$ has only ordinary orbits in $D^\times$.
\end{prop}

\noindent
The proposition immediately follows from the more precise
statement:

\begin{prop}\label{prop:remtr}
Suppose that (\ref{eq:rem}) is verified.
Then the following propositions hold:
\par\noindent
1) There exists a polynomial $P(x)\in K[x]$, with $P(0)\neq 0$, such that the
$q$-difference system $Y(qx)=A_{[P{\mathbb I}_\mu]}(x)Y(x)$
has only trivial singularities in $D^\times$.
\par\noindent
2) Suppose that $({\mathcal S})_q$ has only trivial singularities in $D^\times$.
Then there exists $H(x)\in Gl_\mu(K(x))$ such that
\par\noindent
- $H(x)$ does not have a pole at $0$ and $H(0)\in Gl_\mu(K)$,
\par\noindent
- $H(0)\in Gl_\mu(K)$,
\par\noindent
- $\l\|H(x)\r\|_{0,\rho}^{-1}=\l\|H(x)^{-1}\r\|_{0,\rho}=\rho$,
\par\noindent
- the $q$-difference system $Y(qx)=A_{[H]}(x)Y(x)$ has
only ordinary orbits in $D^\times$.
\end{prop}

\begin{proof}
Let $P(x)\in K[x]$ be a polynomial such that for any
$\zeta\in D^\times$ and any solution matrix $U_\zeta(x)$ meromorphic on a
convenient $q$-invariant analytic domain of center $\zeta$,
the matrix $P(x)U_\zeta(x)$ is analytic at $\zeta$.
Then $P(x)U_\zeta(x)$ is a solution matrix of the $q$-difference system
associated to $A_{[P{\mathbb I}_\mu]}(x)$, which has only trivial singularities in $D^\times$.
This completes the proof of the first part of the statement.
\par
Now we prove {\sl 2)}.
Let $\zeta\in D^\times$ and let
$Y(x)\in M_{\mu\times\mu}({\mathcal A}_{D^\prime})
\cap Gl_\mu({\mathcal M}_{D^\p})$ be a solution at $\zeta$ of $({\mathcal S})_q$.
If $({\mathcal S})_q$ has a trivial singularity at $q^\Z\zeta$, then necessarily
$\det Y(x)$ has a zero in $q^\Z\zeta$. Notice that, since any
infinite subset of $q^\Z\zeta$ has a limit point in $D^\p$,
the analytic function $\det Y(x)\in {\mathcal A}_{D^\prime}$ has
only a finite number of zeros in $q^\Z\zeta$.
\par
Let $q^{n}\zeta$ be a zero of order $k>0$ of $\det Y(x)$
and let $\vec B={}^t(B_1,\dots,B_\mu)$ be a non-zero vector in $K^\mu$ such that $\vec BY(q^n\zeta)=0$.
We fix $1\leq\iota\leq\mu$ such that $|B_\iota|=\max_{j=1,\dots,\mu}|B_j|>0$.
Of course one can suppose that $B_\iota=1$.
Let us consider the gauge transformation matrix
$$
H(x)=\l(\begin{array}{c|c|c}
          {\mathbb I}_{i-1}&\begin{matrix}0\\ \vdots\\ 0\end{matrix} & 0\\
          \hline
          &&\\
          \begin{matrix}\displaystyle\frac{B_1}{x-q^n\zeta}&\cdots&\displaystyle\frac{B_{\iota-1}}{x-q^n\zeta}\end{matrix}
               &\ds\frac{1}{x-q^n\zeta}
               &\begin{matrix}\displaystyle\frac{B_{\iota+1}}{x-q^n\zeta}&
                 \cdots&\displaystyle\frac{B_{\mu}}{x-q^n\zeta}\end{matrix}\\
          &&\\
          \hline
          0&\begin{matrix}0\\ \vdots\\ 0\end{matrix} &{\mathbb I}_{\mu-i}\\
\end{array}\r)\ .
$$
Since $\vec BY(q^n\zeta)=0$, the matrix $H(x)Y(x)$ is still analytic at $q^n\xi$ and
$$
H(x)^{-1}=\l(\begin{array}{c|c|c}
          {\mathbb I}_{i-1}&\begin{matrix}0\\ \vdots\\ 0\end{matrix} & 0\\
          \hline
          &&\\
          \begin{matrix}-B_1&\cdots&-B_{\iota-1}\end{matrix}&x-q^n\zeta
               &\begin{matrix}-B_{\iota+1}&\cdots&-B_{\mu}\end{matrix}\\
          &&\\
          \hline
          0&\begin{matrix}0\\ \vdots\\ 0\end{matrix} &{\mathbb I}_{\mu-i}\\
\end{array}\r)\ .
$$
Moreover $\det \l(H(x)Y(x)\r)$ has a zero at $q^n\zeta$ of order $k-1$.
By iteration, one can construct a basis change satisfying
all the conditions in {\sl 2)}.
\end{proof}
\end{parag}

\subsection{Generic points}
\label{subsec:genericpoints}

We consider an extension $\Omega/K$ of ultrametric fields with the following properties
(for the construction of such a field see for instance
\cite[\S3, 2]{ROB}):

\begin{enumerate}
\setlength{\itemsep}{-2pt}
  \item the field $\Omega$ is complete and algebraically closed;
  \item the set of values of $\Omega$ is $\R_{\geq 0}$;
  \item the residue field of $\Omega$ is transcendental over the
residue field of $K$;
  \item for any $R\in\R_{\geq 0}$ the field $\Omega$ contains an element $t_R$,
transcendent over $K$,  such that
$|t_R|=R$ and that the norm induced by $\Omega$ over $K(t_R)$ is defined by
$$
\l|\frac{\sum a_i t_R^i}{\sum b_j t_R^j}\r|=
\frac{\sup_i |a_i|R^i}{\sup_j |b_j|R^j}\ .
$$
\end{enumerate}
Observe that, if $f(x)$ is an analytic function over a disk of center $0$ and
radius $\rho>R$, we have
$$
|f(t_R)|=\l\|f(x)\r\|_0(R)=\sup_{x\in\Omega, |x|< R}|f(x)|\ .
$$

\begin{defn}
We call $t_R$ a \emph{generic point with respect to $K$
(at distance $R$ from $0$)}.
A \emph{generic point $t_{\xi,R}$ with respect to $K$ (at distance $R$ from $\xi\in K$)}
is defined by shifting.
\end{defn}

\medskip
Dwork-Robba's generic points play a fundamental role in $p$-adic differential equation theory:
historically, their introduction has been the first attempt to fill the gap left by the
absence of a $p$-adic analytic continuation.
In fact, the radius of solutions at $t_{\xi,R}$ is a sort of global invariant for the differential equation
that is equal to the radius of convergence of solutions at almost any point
of the disk $D(\xi,R^-)$ and allows for an estimate in the other points.
From a more recent point of view, one should think of generic points as points
of a Berkovich analytic space.

\medskip
Consider a $q$-difference system of order $\mu$ with meromorphic coefficients over a
$q$-invariant disk $D(\xi,\rho^-)$:
$$
Y(qx)=A(x)Y(x)\ ,\
\hbox{with }A(x)\in Gl_\mu({\mathcal M}_{D,K}).
\leqno{(\mathcal S)_q}$$
As in the previous section (cf. (\ref{eq:analyticsolution})),
for any $R<\rho$ one can consider the formal solution
of $(\mathcal S)_q$ at $t_{\xi,R}$:
\begin{equation}\label{eq:solgen}
Y(t_{\xi,R},x)=\sum_{n\geq 0}\frac{G_n(t_{\xi,R})}{[n]_q^!}(x-t_{\xi,R})_{q,n}
\in Gl_\mu(\Omega\[[x-t_{\xi,R}\]]_q)\ .
\end{equation}
Notice that since the system $(\mathcal S)_q$ is defined over $K$,
the matrix $A(x)$ cannot have any poles in $q^\N t_{\xi,R}$ and
$\det A(x)$ cannot have any zeros in $q^\N\xi$: it follows that
$Y(t_{\xi,R},x)$ is necessarily in $Gl_\mu(\Omega\[[x-t_{\xi,R}\]]_q)$.

\begin{defn}
We call \emph{generic radius of convergence of $({\mathcal S})_q$ at $t_{\xi,R}$}
the number
$$
\chi_{\xi,R}(A,q)=\inf\l(R,
\liminf_{n\rightarrow\infty}\l|\frac{G_n(t_{\xi,R})}{[n]_q^!}\r|^{-1/n}\r)\ .
$$
\end{defn}

\begin{lemma}\label{lemma:genradius}
Let $1\geq R\geq |(1-q)\xi|$.
The generic radius of convergence $\chi_{\xi,R}(A,q)$ is
invariant under meromorphic gauge transformation and
$$
\chi_{\xi,R}(A,q)\geq\frac{R|\pi_q|}{\sup\l(1,|G_1(t_{\xi,R})|\r)}\ .
$$
\end{lemma}

\begin{proof}
The first assertion is proved in \cite[4.2.3]{DV} in the case $\xi=0$ and $R=1$,
but the same proof applies to this case. The second assertion is immediately deduced from the recursive
relation satisfied by the $G_n(x)$'s and (\ref{eq:normaqder}), since $\lim_{n\rightarrow\infty}|[n]_q^!|^{1/n}=|\pi_q|$.
\end{proof}

\begin{lemma}\label{lemma:transferevidente}
If $A(x)$ is analytic over the disk
$D(\xi,\rho^-)$ and $\rho\geq\chi_{\xi,R}(A,q)>|(1-q)\xi|$, then
$({\mathcal S})_q$ has a fundamental analytic solution over $D(\xi,\chi_{\xi,R}(A,q)^-)$.
\end{lemma}

\begin{proof}
It is enough to remark that
$$
\l|\frac{G_n(\xi)}{[n]_q^!}\r|\leq
\l|\frac{G_n(t_{\xi,R})}{[n]_q^!}\r|\ .
$$
\end{proof}

The theorems estimating the radius of convergence of $Y(\xi,x)$ with respect to
$\chi_{\xi,R}(A,q)$ are usually called {\it transfer theorems}:
in the next chapter we will prove a transfer theorem from a disk where an
analytic solution exists to a contiguous disk, where the system has only ordinary orbits.
This result is a consequence of the effective bound theorem.
Chapters \ref{sec:frobenius} and \ref{sec:transfer} are devoted to the proof of a transfer theorem for regular singular disks.

\begin{parag}\label{parag:dworkfrobenius}
{\bf The cyclic vector lemma and the $q$-analogue of the Dwork-Frobenius theorem.}

\par\medskip\noindent
It may seem that calculating a generic radius of convergence is as difficult as
calculating a radius of convergence at points which are rational over $K$.
This is not completely true; in fact the generic radius of convergence
is very easy to calculate when it is small, using a
$q$-difference equation associated to $({\mathcal S})_q$: it is the $q$-analogue of the Dwork-Frobenius theorem
\cite[VI, 2.1]{DGS}.

\par\medskip
As in the differential world, a $q$-difference equation associated to $Y(qx)=A(x)Y(x)$,
with $A(x)\in Gl_\mu({\mathcal M}_{D,K})$, is constructed using a cyclic vector lemma
(cf. for instance \cite[Annexe B]{Sfourier} or \cite[1.3]{DV}), which states the existence of a meromorphic matrix
$H(x)\in Gl_\mu({\mathcal M}_{D,K})$ such that
\begin{equation}\label{eq:cyclicvector}
A_{[H]}=
\l(\begin{array}{c|c}
0& \\
\vdots &{\mathbb I}_{\mu-1}\\
0& \\
\hline\\
a_0(x)&a_1(x)\dots a_{\mu-1}(x)
\end{array}\r)\ .
\end{equation}
Then $y(x)$ is a solution of the $q$-difference equations
\begin{equation}\label{eq:equation}
y(q^\mu x)-a_{\mu-1}(x)y(q^{\mu-1}x)-\dots-a_0(x)y(x)=0
\end{equation}
if and only if
$$
\begin{pmatrix}y(qx)\\y(q^2 x)\\ \vdots\\ y(q^\mu x)\end{pmatrix}
=A_{[H]}\begin{pmatrix}y(x)\\y(qx)\\ \vdots\\ y(q^{\mu-1}x)\end{pmatrix}\ .
$$
Consider the lower triangular gauge transformation matrix
\begin{equation}\label{eq:cyclicgauge}
\widetilde H=\l(a_{i,j}\r)_{i,j=0,\dots,\mu-1}\ ,
\hbox{ with }
a_{i,j}=\begin{cases}
\ds\frac{1}{x^i}\frac{(-1)^{i+j}}{(q-1)^i}{i\choose j}_{q^{-1}}q^{-\frac{j(j-1)}{2}}& \hbox{ if }j\leq i;\\ \\
0&\hbox{ otherwise}.
\end{cases}
\end{equation}
It follows from (\ref{eq:leibnizquattro}) that
$$
G_{[\tilde HH]}=\frac{A_{[\tilde HH]}-{\mathbb I}_\mu}{(q-1)x}
=\l(\begin{array}{c|c}
0& \\
\vdots &{\mathbb I}_{\mu-1}\\
0& \\
\hline\\
b_0(x)&b_1(x)\dots b_{\mu-1}(x)
\end{array}\r)\
$$
and that $y(x)$ is a solution of the $q$-difference equations (\ref{eq:equation})
if and only if
$$
\begin{pmatrix}\dq y(x)\\\dq^2 y(x)\\ \vdots\\\dq^\mu y(x)\end{pmatrix}
=G_{[\tilde HH]}\begin{pmatrix}y(x)\\\dq y(x)\\ \vdots\\\dq^{\mu-1}y(x)\end{pmatrix}\ .
$$

\begin{prop}\label{prop:dworkfrobenius}
Let  $|1-q|<|\pi|$.
If $\sup_{i=0,\dots,\mu-1}|b_i(t_{0,R})|>R^{i-\mu}$ then
$$
\chi_{0,R}(A,q)=\frac{|\pi|}
{\sup_{i=0,\dots,\mu-1}|b_i(t_{0,R})|^{1/(\mu-i)}}\ .
$$
\end{prop}

\begin{proof}
The proposition is proved in \cite[4.3]{DV} in the case $R=1$, $\xi=0$.
We deduce the statement above by rescaling.
\end{proof}

\end{parag}

\makeatletter\newcount\uffa\uffa=\c@subsection\makeatother
\section{Effective bounds for $q$-difference systems}
\label{sec:ee}
\makeatletter\c@subsection=\uffa\makeatother

In this chapter we prove an effective bound theorem for $q$-difference
systems: it is the analogue of a theorem of Dwork and Robba
(cf. \cite{DR} for the proof in the case of differential equations.
The statement concerning differential systems is proved for instance in
\cite{B}, \cite{AGfunctions} and \cite{DGS}).
\par
Let us explain the effective bound theorem for an analytic differential
equations of order one:
the theorem is actually almost trivial in this case, but
we can already point out the differences with the $q$-difference version.
\par
Let $\xi\in{\mathbb A}^1_K$, $\xi\neq 0$, and let $u(x)$ be a meromorphic
function over an open disk of center $\xi$ and radius $\rho>0$.
For all $R\in (0,\rho)$, the multiplicative norm
$\|~\|_\xi(R)$ (cf. \S\ref{subsec:gaussnorm})
induces a norm over the field of meromorphic
functions over $D(\xi,\rho^-)$.
Let
$$
g_n(x)=\frac{1}{n!}\frac{d^nu}{dx^n}(x)u(x)^{-1}\ ,\ \forall n\geq 0\ .
$$
The effective bound theorem for differential equation of order $1$ states that
for any $R\in (0,\rho)$
$$
\|g_n(x)\|_\xi(R)\leq R^{-n}\ .
$$
Of course this inequality is easy to prove, in fact
\begin{equation*}
  \begin{array}{rcl}
         \|g_n(x)\|_\xi(R)
         &\leq&\ds\l\|\frac{1}{n!}\frac{d^nu}{dx^n}(x)\r\|_\xi(R)\l\|u^{-1}(x)\r\|_\xi(R)\\
         &\leq&\ds R^{-n}\|u(x)\|_\xi(R)\|u(x)\|_\xi(R)^{-1}\\
         &\leq&\ds R^{-n}\ .
  \end{array}
\end{equation*}
The multiplicativity of $\|~\|_\xi(R)$  is the key point of the inequalities
above.
\par
Let us consider a me\-ro\-morphic $q$-difference system: its solutions may be
meromorphic over a non-connected analytic domain
(as the one considered in (\ref{eq:primoqdisco})).
Therefore a natural analogue
of $\|~\|_\xi(R)$ would be a sup-norm over a non-connected domain,
which is necessarily non-multiplicative.
In particular $\|u(x)^{-1}\|_\xi(R)$ could be greater than $\|u(x)\|_\xi(R)^{-1}$.
In other words, the assumption of having a solution matrix (analytic or meromorphic) over an
$q$-invariant disk cannot be avoided.

\subsection{Effective bound theorem for $q$-difference systems}

Let $D=D(\xi,\rho^-)$ be an open disk of center $\xi\in{\mathbb A}^1_K$ and
${\mathcal M}_D$ be the field of meromorphic functions over $D$.
For any $R\in(0,\rho)$ the norm $\|~\|_\xi(R)$ extends from ${\mathcal A}_D$ to
${\mathcal M}_D$ by multiplicativity. For any $f(x)\in {\mathcal M}_D$ one usually sets:
$$
\|f(x)\|_{\xi,\rho}=\lim_{R\rightarrow\rho}\|f(x)\|_\xi(R)\ .
$$
This limit may be not bounded, but respects multiplication of functions, as
well as $\|~\|_\xi(R)$.

\par
We suppose that $D$ is $q$-invariant and we consider $Y(x)\in Gl_\mu({\mathcal M}_D)$.
Then $Y(x)$ is solution of the $q$-difference system
defined by
$$
G_n(x)=\dq^nY(x)Y(x)^{-1}\ ,\ \forall n\geq 0.
$$
Obviously the entries of $G_n(x)$ are meromorphic functions over $D$.

\begin{thm}\label{thm:eesys}
$$
\l\|\frac{G_n(x)}{[n]_q^!}\r\|_{\xi,\rho}\leq \{n,\mu-1\}_p^q
\l(\sup_{i=0,\dots,\mu-1}\l\|G_i(x)\r\|_{\xi,\rho}\rho^i\r)\frac{1}{\rho^n}\ ,
$$
where
\begin{equation*}
  \{n,\mu-1\}_p^q=
  \begin{cases}
      1, &\hbox{if n=0,}\\
      \ds\sup_{1\leq\la_1<\dots<\la_{\mu-1}\leq n, \la_i\in\Z} \frac{1}{\l|[\la_1]_q\cdots[\la_{\mu-1}]_q\r|},
      &\hbox{otherwise.}
  \end{cases}
\end{equation*}
\end{thm}

One can also consider a matrix $Y(x)$ analytic over a $q$-invariant disk
around $t_{\xi,R}$, with coefficients in $\Omega$ (cf. \S\ref{subsec:genericpoints}),
and assume that it is solution of a $q$-difference system
$({\mathcal S})_q$ defined over $K$. Then:

\begin{cor}
$$
\l|\frac{G_n(t_{\xi,R})}{[n]_q^!}\r|
\leq \{n,\mu-1\}_p^q \l(\sup_{i=0,\dots,\mu-1}\l|G_i(t_{\xi,R})\r|
\chi_{\xi,R}(A,q)^i\r)
\frac{1}{\chi_{\xi,R}(A,q)^n}\ .
$$
\end{cor}

\medskip
The proof of theorem \ref{thm:eesys} (cf. (\ref{parag:eesys}) below)
follows the proof of the Dwork-Robba theorem concerning
the effective bounds for $p$-adic differential systems. As in the
differential case, it relies on the analogous result
for $q$-difference equations.

\begin{parag}
{\bf Effective bounds for $q$-difference equations.}
\label{parag:eeequa}

\par\medskip\noindent
We consider $\vec u=(u_1,\dots,u_\mu)\in{\mathcal M}_D^\mu$, with $\mu\in\Z$,
$\mu\geq 1$, such that $u_1,\dots,u_\mu$
are linearly independent over the field of constants $K$ of ${\mathcal M}_D$.
By the $q$-analogue of the Wronskian lemma
(cf. \cite[Appendice]{Sfourier} or \cite[\S1.2]{DV}) this is equivalent to
supposing that the $q$-Wronskian matrix
\begin{equation*}
W_q(\vec u)=\begin{pmatrix}
                  u_1(x) & \dots & u_\mu(x) \\
                  \dq u_1(x) & \dots & \dq u_\mu(x) \\
                  \vdots & \ddots &\vdots \\
                  \dq^{\mu-1} u_1(x) & \dots & \dq^{\mu-1} u_\mu(x)
             \end{pmatrix}
\end{equation*}
is in $Gl_\mu({\mathcal M}_D)$.
\par
We consider the $q$-difference equation defined by
$$
\dq^n\vec u=\vec g_n(x)W_q(\vec u)\ ,
\hbox{ where $\vec g_n\in{\mathcal M}_D^\mu$ and $n\geq 0$.}
$$

\begin{lemma}\label{lemma:eeequa}
  \begin{equation}\label{eq:eeequa}
    \l\|\frac{\vec g_n}{[n]_q^!}\r\|_{\xi,\rho}\leq \{n,\mu-1\}_p^q\frac{1}{\rho^n}\ .
  \end{equation}
\end{lemma}

\begin{proof}
By extending $K$, we can find  $\a\in K$ such that $|\a|=\rho$.
If we set $y=\a x$ and $\zeta=\a \xi$, we are reduced to proving the lemma
for $\rho=1$.
\par
So we suppose $\rho=1$ and we prove the lemma by induction on $\mu$.
For $\mu=1$ and any $R\in((q-1)\xi,1)$ we have
(cf. (\ref{eq:normaqder}))
$$
\l\|\frac{\vec g_n}{[n]_q^!}\r\|_\xi(R)
\leq\l\|\frac{\dq^n u_1}{[n]_q^!}(x)\r\|_\xi(R)\|u_1(x)^{-1}\|_\xi(R)
\leq \frac{1}{R^n}\ .
$$
Letting $R\rightarrow 1$, we obtain (\ref{eq:eeequa}).
\par
Let $\mu>1$. For $n=0$ the inequality is trivial, so let $n>0$.
We set:
$$
\vec u=(u_1,\dots,u_{\mu+1})=u(1,\vec\tau)\in{\mathcal M}_D^{\mu+1}\ ,
\hbox{ where $u=u_1$ and $\vec\tau=(\tau_1,\dots,\tau_\mu)\in{\mathcal M}_D^\mu$.}
$$
The idea of the proof is to apply the inductive hypothesis to $\dq\vec\tau$.
We know that the vector $\vec h_n\in{\mathcal M}_D^\mu$,
defined for any $n\geq 0$ by
$$
\dq^n\l(\dq\vec\tau\r)=\vec h_n(x)W_q\l(\dq\vec\tau\r)\ ,
$$
satisfies the inequality
\begin{equation}\label{eq:bleark1}
  \l\|\frac{\vec h_n}{[n]_q^!}\r\|_{\xi,1}\leq\{n,\mu-1\}_p^q\ .
\end{equation}
Notice that $\vec h_n$ and $\vec g_n$ verify the relation
\begin{equation}\label{eq:bleah}
  \begin{array}{rcl}
    \lefteqn{\ds\frac{\vec g_n}{[n]_q^!}W_q(\vec u)
    =\ds\frac{\dq^n}{[n]_q^!}\vec u=\frac{\dq^n}{[n]_q^!}\l(u(1,\vec\tau)\r)}\\
    &=&\ds\frac{\dq^n}{[n]_q^!}(u)(1,\vec\tau)+
    \l(0,\sum_{j=1}^n\frac{\dq^{n-j}}{[n-j]_q^!}(u)(q^jx)\frac{1}{[j]_q}
    \frac{\vec h_{j-1}}{[j-1]_q^!}\Big(0, W_q(\dq\vec\tau)\Big)\r)\ ,
  \end{array}
\end{equation}
while $W_q(\vec u)$ and $W_q(\dq\vec\tau)$ satisfy
$$
W_q(\vec u)=W_q\l(u(1,\vec\tau)\r)=uP\begin{pmatrix}1 &\vec\tau\\  0&
  W_q(\dq\vec\tau)\end{pmatrix}\ ,
$$
where
$$
P=\frac{1}{u(x)}
\begin{pmatrix}
u(x)    &0
    &0
    &\cdots
    &0     \\
\dq u(x)&u(qx)
    &0
    &\ddots
    &\vdots\\
{2\choose 2}_q\dq^2 u(x)  &{2\choose 1}_q\dq u(qx)
    &u(q^2x)
    &\ddots
    &0\\
\vdots  &\vdots
    &\vdots
    &\ddots
    &0\\
{\mu\choose\mu}_q\dq^\mu u(x)  &{\mu\choose\mu-1}_q\dq^{\mu-1}u(qx)
    &{\mu\choose\mu-2}_q\dq^{\mu-2}u(q^2x)
    &\cdots
    &u(q^\mu x)
\end{pmatrix}\ .
$$
Since $\l\|{\mu\choose\mu-i}_q\dq^{\mu-i}u(q^ix)u(x)^{-1}\r\|_{\xi,1}\leq 1$
and $\l\|u(q^i x)u(x)^{-1}\r\|_{\xi,1}=1$, we have
$\|P\|_{\xi,1}=\|\det P\|_{\xi,1}=1$ and hence $\|P^{-1}\|_{\xi,1}\leq 1$.
This implies that
\begin{equation}\label{eq:bleark2}
\l\|
     \begin{pmatrix}
       1 &\vec\tau\\
       0& W_q(\dq\vec\tau)
     \end{pmatrix}
W_q(\vec u)^{-1}
\r\|_{\xi,1}
=\|u^{-1}P^{-1}\|_{\xi,1}
\leq\|u^{-1}\|_{\xi, 1}\ .
\end{equation}
We obtain the desired inequality by combining (\ref{eq:bleark1}), (\ref{eq:bleah}) and (\ref{eq:bleark2}):
\begin{equation*}
  \begin{array}{rcl}
    \ds\l\|\frac{\vec g_n}{[n]_q^!}\r\|_{\xi,1}
    &\leq&\ds\sup_{j=1,\dots,n}\l(1,\l\|\frac{\vec h_{j-1}}{[j]_q^!}\r\|_{\xi,1}\r)\\
    &\leq&\ds\sup_{j=1,\dots,n}\l(1,\l\|\frac{1}{[j]_q}\r\|_{\xi,1}\{j-1,\mu-1\}_p^q\r)\\ \\
    &\leq&\ds\{j,\mu-1\}_p^q\ .
  \end{array}
\end{equation*}
\end{proof}
\end{parag}

\begin{parag}
{\bf Proof of theorem \ref{thm:eesys}.}
\label{parag:eesys}

\par\medskip\noindent
We recall that we are given a matrix $Y(x)\in Gl_\mu({\mathcal M}_D)$ meromorphic
over a $q$-invariant open disk $D$ and that we
have set
$$
G_n(x)=\dq^nY(x)Y(x)^{-1}\in M_{\mu\times\mu}({\mathcal M}_D)\ .
$$
We want to prove that
$$
\l\|\frac{G_n(x)}{[n]_q^!}\r\|_{\xi,\rho}\leq \{n,\mu-1\}_p^q
\l(\sup_{i=0,\dots,\mu-1}\l\|G_i(x)\r\|_{\xi,\rho}\rho^i\r)\frac{1}{\rho^n}\ .
$$
First of all we observe that:
\par\noindent
1) as in the case of $q$-difference equations it is enough to prove the
inequality for $\rho=1$;
\par\noindent
2) it is enough to prove the inequality above for the first row of $G_n(x)$.
\par\noindent
So we suppose $\rho=1$ and we call $\vec u=(u_1,\dots,u_\mu)\in{\mathcal M}_D^\mu$ the
first row of $Y(x)$. Let $k\leq \mu$ be the rank of $W_q(\vec u)$
and let $E\in Gl_\mu(K)$ be such that
$$
\vec uE=(\vec z,\underline 0)\ ,\hbox{ with }\vec z\in{\mathcal M}_D^k
\hbox{ and }\underline 0=(\underbrace{0,\dots,0}_{\times (n-k)})\ .
$$
By (\ref{lemma:eeequa}) the vectors $\vec h_n$, defined by
$$
\dq^n(\vec z)=\vec h_n W_q(\vec z)\ ,
\hbox{ for any $n\in\Z_{\geq 0}$,}
$$
satisfy the inequality
$$
\l|\frac{\vec h_n}{[n]_q^!}\r|_{\xi,1}\leq \{n,\mu-1\}_p^q\ .
$$
Moreover the $q$-wronskian $W_q(\vec z)$ satisfies
$$
\begin{pmatrix}\vec u\\ \dq\vec u\\ \vdots\\ \dq^{k-1}\vec u\end{pmatrix}E
=\l(W_q(\vec z),\underline 0\r)\in M_{k\times\mu}({\mathcal M}_D)\ .
$$
In order to deduce an estimate of $\frac{G_n(x)}{[n]_q^!}$ from the
above estimate of $\frac{\vec h_n}{[n]_q^!}$ we need to
express $\frac{\dq^n}{[n]_q^!}(\vec u)$ in terms of $\frac{\vec h_n}{[n]_q^!}$:
$$
\frac{\dq^n}{[n]_q^!}(\vec u) E
=\frac{\dq^n}{[n]_q^!}\l(\vec u E\r)=\frac{\dq^n}{[n]_q^!}\l(\vec z,\underline 0\r)
    =\frac{\vec h_n}{[n]_q^!}\l(W_q(\vec z),\underline 0\r)
=\frac{\vec h_n}{[n]_q^!}
\begin{pmatrix}\vec u\\ \dq\vec u\\ \vdots\\ \dq^{k-1}\vec u\end{pmatrix}E\ .
$$
Since $E\in Gl_\mu(K)$, we obtain
$$
\l(\hbox{first row of $\ds\frac{G_n(x)}{[n]_q^!}$}\r)Y(x)
=\frac{\dq^n}{[n]_q^!}(\vec u)
=\frac{\vec h_n}{[n]_q^!}
\begin{pmatrix}\vec u\\ \dq\vec u\\ \vdots\\ \dq^{k-1}\vec u\end{pmatrix}\
$$
and hence we deduce that
\begin{equation*}
  \begin{array}{rl}
    \l\|\hbox{first row of $\ds\frac{G_n(x)}{[n]_q^!}$}\r\|_{\xi,1}
    &\leq\l\|\frac{\vec h_n}{[n]_q^!}\r\|_{\xi,1}
    \l\|\begin{pmatrix}\vec u\\ \dq\vec u\\ \vdots\\
    \dq^{k-1}\vec u\end{pmatrix}Y(x)^{-1}\r\|_{\xi,1}\\
    &\ds\leq\{n,\mu-1\}_p^q\l(\sup_{i=0,\dots,k-1}\l\|\hbox{first row of
      $G_i(x)$}\r\|_{\xi,1}\r)\ .
  \end{array}
\end{equation*}
This completes the proof.
\hfill$\square$
\end{parag}

\subsection{Some consequences: a transfer theorem in
ordinary disks and a corollary about $q$-deformations}
\label{subsec:consequences}

\begin{parag}
{\bf Transfer theorem in ordinary disks}

\par\medskip\noindent
As in $p$-adic differential equation theory,
the following transfer result follows from the effective bound
theorem:
\end{parag}

\begin{cor}\label{cor:transferfacile}
Let $G(x)\in M_{\mu\times\mu}({\mathcal M}_{D(\xi,\rho^-)})$.
We suppose that the $q$-difference system
$\dq Y(x)=G(x)Y(x)$ has a meromorphic fundamental solution
over the $q$-invariant disk
$D(\xi,\eta^-)$, with $\eta<\rho$, and that there exists
a point $\zeta$ such that
\par\noindent
- $|\xi-\zeta|=\eta$,
\par\noindent
- $G(x)$ is analytic over $D(\zeta,\eta^-)$.
\par
Then $\dq Y(x)=G(x)Y(x)$ has an analytic solution in $D(\zeta,\eta^-)$.
\end{cor}

\begin{proof}
Let $\dq^n  Y(x)=G_n(x)Y(x)$ for any non-negative integer $n$.
By the previous theorem we have
$$
\l|\frac{G_n(t_{\xi,\eta})}{[n]_q^!}\r|=\l\|\frac{G_n(x)}{[n]_q^!}\r\|_{\xi,\eta}
\leq C\{n,\mu-1\}_p^q\frac{1}{\eta^n}\ ,
$$
where $C$ is a constant depending only on $G(x)$ and $\eta$.
The point $t_{\xi,\eta}$ is also a generic point at distance $\eta$ from $\zeta$, hence
$$
\l|\frac{G_n(\zeta)}{[n]_q^!}\r|\leq\l|\frac{G_n(t_{\xi,\eta})}{[n]_q^!}\r|\ .
$$
Let $\kappa$ be the smallest positive integer such that $|1-q^\kappa|<|\pi|$
and let $p^{l_n}$ be the greatest integer power of $p$ smaller or equal to $n$.
Then we have
$$
\{n,\mu-1\}_p^q\leq \l|\l[p^{l_n}\kappa\r]_q\r|^{1-\mu}
\leq\l(|p^{l_n}[\kappa]_q|\r)^{1-\mu}
=\l(p^{-l_n}|[\kappa]_q|\r)^{1-\mu}
\leq n^{\mu-1}|[\kappa]_q|^{1-\mu}\ .
$$
We conclude by applying (\ref{prop:taylor}), since
$$
\limsup_{n\rightarrow\infty}\l|\frac{G_n(\zeta)}{[n]_q^!}\r|^{1/n}
\leq \eta^{-1}
\limsup_{n\rightarrow\infty}n^{(\mu-1)/n}|[\kappa]_q|^{-(\mu-1)/n}=\eta^{-1}\ .
$$
\end{proof}

\begin{parag}\label{parag:deformation}
{\bf Effective bounds and $q$-deformation of $p$-adic differential equations}

\par\medskip\noindent
Let $q_k\in K$, $k\in\N$, be a sequence such that $q_k\rightarrow 1$ when
$k\rightarrow \infty$ and
$G^{(k)}(x)$ a sequence of square matrix of order $\mu$ whose entries are
\emph{analytic bounded functions} over an open disk $D$ of center $\xi$ and radius $\rho$,
with $|\xi|\leq\rho\leq 1$.
By rescaling, we can assume that $D=D(\xi,1^-)$ and $|\xi|\leq 1$.
Suppose that the sequence of matrices $G^{(k)}(x)$ tends to a matrix $G(x)$ uniformly over
any $D(\xi,\eta^+)\subset D$.
Then we say that the family of systems
$$
d_{q_k}Y(x)=G^{(k)}(x)Y(x)
\leqno{({\mathcal S})^\p_{q_k}}
$$
is a \emph{$q$-deformation of the differential system}
$$
\frac{dY}{dx}(x)=G(x)Y(x)\ .
\leqno{({\mathcal S})}$$
Let $G^{(k)}_n(x)$ and $G_n(x)$, for $n\geq 0$, be the matrices respectively defined by
$$
d_{q_k}^nY(x)=G^{(k)}_n(x)Y(x)
\hbox{ and }
\frac{d^nY}{dx^n}(x)=G_n(x)Y(x)\ .
$$
We assume that
$$
\|G(x)\|_{\xi,1}\leq 1\ .
$$
\end{parag}

\begin{lemma}
Let $\veps\in\R_{>0}$.
There exists $k_\veps>>0$ such that,
if
$$
\|G(x)-G^{(k)}(x)\|_{\xi,1}<\veps,
\hbox{ for any $k\geq k_\veps$, }
$$
then
$$
\|G_n(x)-G^{(k)}_n(x)\|_{\xi,1}<\veps,
\hbox{ for any $k\geq k_\veps$ and any $n\geq 1$. }
$$
\end{lemma}

\begin{proof}
We want to prove the statement by induction on $n$, the case $n=1$ being true by assumption.
\par
Let $k_\veps$ be a positive integer such that $|1-q_k|<\veps$
for all $k\geq k_\veps$.
Suppose
$$
\|G^{(k)}_n(x)-G_n(x)\|_{\xi,1}<\veps\ .
$$
Then by  (\ref{lemma:stimaqdef}), for any $k\geq k_\veps$, we have
$$
\l\|d_{q_k}G^{(k)}_n(x)-\frac{d}{dx}G_n(x)\r\|_{\xi,1}<\veps\ .
$$
Moreover $G^{(k)}_n(x)$ are analytic bounded function over $D(\xi, 1^-)$,
hence
$$
\l\|G^{(k)}_n(q_kx)-G^{(k)}_n(x)\r\|_{\xi,1}\leq |1-q_k|<\veps\ .
$$
Finally
$$
\begin{array}{rcl}
\lefteqn{\|G^{(k)}_{n+1}(x)-G_{n+1}(x)\|_{\xi,1}}\\ \\
&=   &\l\|G^{(k)}_1(x)G^{(k)}_n(q_kx)-d_{q_k}G^{(k)}_n(x)-G_1(x)G_n(x)+\frac{d}{dx}G_n(x)\r\|_{\xi,1}\\ \\
&\leq&\sup\Big(\l\|(G^{(k)}_1(x)-G_1(x))G^{(k)}_n(q_kx)\r\|_{\xi,1},
    \l\|G_1(x)(G^{(k)}_n(q_kx)-G^{(k)}_n(x))\r\|_{\xi,1},\\ \\
&    &\l\|G_1(x)(G^{(k)}_n(q_kx)-G_n(x))\r\|_{\xi,1},
    \l\|d_{q_k}G^{(k)}_n(x)-\frac{d}{dx}G_n(x)\r\|_{\xi,1}\Big)\\ \\
&<   &\veps\ .
\end{array}
$$
\end{proof}

\begin{prop}\label{prop:qdeformazione}
Under the assumption above we have:
\par\noindent
1) For any $k>>0$ the $q_k$-difference system $d_{q_k}Y(x)=G^{(k)}(x)Y(x)$
has an analytic fundamental solution $Y^{(k)}(x)$ over a disk $D(\xi,\eta_k^-)\subset D$,
verifying $Y^{(k)}(\xi)={\mathbb I}_\mu$. Moreover, $\liminf_{k\rightarrow\infty}\eta_k>0$.
\par\noindent
2) Let $\eta=\liminf_{k\rightarrow\infty}\eta_k$. Then $Y^{(k)}(x)$ tends pointwise over $D(\xi,\eta^-)$
to a fundamental solution $Y(x)$ of $\frac{dY(x)}{dx}=G(x)Y(x)$.
Moreover $Y^{(k)}(x)$ tends uniformly to $Y(x)$ over $D(\xi,(\eta^\p)^-)$, for any $0<\eta^\p<\eta$.
\end{prop}

\begin{proof}\hfill\\
\smallskip\noindent
{\sl 1)}
By (\ref{lemma:solformali}) the system $({\mathcal S})^\p_{q_k}$ has a formal solution of the form
$$
Y^{(k)}(\xi,x)=\sum_{n\geq 0}\frac{G^{(k)}_n(\xi)}{[n]_{q_k}^!}(x-\xi)_{q_k,n}\ .
$$
We set
$$
\eta_k^{-1}=\limsup_{n\rightarrow\infty}\l|\frac{G^{(k)}_n(\xi)}{[n]_{q_k}^!}\r|^{1/n}\ .
$$
Since $G^{(k)}(x)$ converges unifomely to $G(x)$ over any closed
disk contained in $D$, there exists a real positive constant $C$
such that
$$
\l\|G_n^{(k)}(x)\r\|_\xi(1)\leq\sup\l(1,\r\|G^{(k)}(x)\l\|_\xi(1)\r)^n
\leq C^n
$$
for any $k>>0$. Therefore for $k>>0$ we have
$\eta_k\geq|\pi_q|C^{-1}\geq|\pi|C^{-1}$.
It follows that $\eta=\liminf_{k\to\infty}\eta_k>0$ and that
for $k>>0$ such that $|(1-q_k)\xi|<|\pi|C^{-1}$ the matrix $Y^{(k)}(\xi,x)$ is the $q$-expansion of an
analytic fundamental solution of $d_{q_k}Y(x)=G^{(k)}(x)Y(x)$.

\par\smallskip\noindent
{\sl 2)} The proof of the second assertion follows faithfully the proof of \cite[IV, 5.4]{DGS}.
Let $\eta_k$ be the radius of convergence of $Y^{(k)}(\xi,x)$.
Letting $k\rightarrow\infty$ in the effective bound estimate
$$
\begin{array}{l}
\ds\l|\frac{G^{(k)}_n(\xi)}{[n]_{q_k}^!}\r|\leq\l\{n,\mu-1\r\}_p^{q_k}C_k\frac{1}{\eta_k^n}\ ,
\ds\hbox{ where }C_k=\sup_{i=0,\dots,\mu-1}\l\|G^{(k)}_i(x)\r\|_{\xi,\eta_k}\eta_k^i\ ,
\end{array}
$$
we obtain (cf. \cite[IV, 3.1]{DGS})
$$
\begin{array}{l}
\ds\l|\frac{G_n(\xi)}{n!}\r|\leq \l\{n,\mu-1\r\}_pC\frac{1}{\eta^n}\ ,
\ds\hbox{ where } C=\sup_{i=0,\dots,\mu-1}\l\|G_i(x)\r\|_{\xi,\eta}\eta^i\ ,\\
\hbox{ and }\{n,\mu-1\}_p=
  \begin{cases}
      1, &\hbox{if n=0,}\\
      \ds\sup_{1\leq\la_1<\dots<\la_{\mu-1}\leq n} \frac{1}{\l|\la_1\cdots\la_{\mu-1}\r|},
      &\hbox{otherwise.}
  \end{cases}
\end{array}
$$
This proves that the solution $Y(\xi,x)$ of $({\mathcal S})_q$ at $\xi$ converges over $D(\xi,\eta^-)$.
\par
As far as the uniform convergence is concerned, notice that for any $R\in(0,\eta)$
and $\veps>0$ there exist $\bar n,\bar k\in\N$ such that
for any $n\geq \bar n$ and any $k\geq \bar k$ we have:
$$
\l\|\frac{G^{(k)}_n(\xi)}{[n]_{q_k}^!}(x-\xi)_{q_k,n}\r\|_\xi(R)
\leq C_k\{n,\mu-1\}_p^{q_k}\l(\frac{R}{\eta_k}\r)^n
\leq C_k n^{\mu-1}\l(\frac{R}{\eta_k}\r)^n<\veps
$$
and
$$
\l\|\frac{G_n(\xi)}{n!}(x-\xi)^n\r\|_\xi(R)\leq C\{n,\mu-1\}_p\l(\frac{R}{\eta}\r)^n
\leq Cn^{\mu-1}\l(\frac{R}{\eta}\r)^n<\veps
$$
By considering a bigger $\bar k$ and using an analogous
argument to the one in the proof of (\ref{lemma:stimaqdef}),
we can suppose that for any $k\geq \bar k$ we have also
$$
\begin{array}{llr}
\lefteqn{\ds\l\|\sum_{n=0}^{\bar n}\frac{G^{(k)}_n(\xi)}{[n]_{q_k}^!}(x-\xi)_{q_k,n}-
    \sum_{n=0}^{\bar n}\frac{G_n(\xi)}{n!}(x-\xi)^n\r\|_\xi(R)}\\
&\leq &\ds\sup_{n=0,\dots,\bar n}\l\|\frac{G^{(k)}_n(\xi)}{[n]_{q_k}^!}(x-\xi)_{q_k,n}-
    \frac{G_n(\xi)}{n!}(x-\xi)^n\r\|_\xi(R)<\veps\ .
\end{array}
$$
Finally we obtain the uniform convergence over any $D(\xi, R^-)$:
$$
\begin{array}{rcl}
\lefteqn{\ds\big\|Y(\xi,x)-Y^{(k)}(\xi,x)\big\|_\xi(R)}\\
&\leq&\ds\sup\l(\l\|\sum_{n=0}^{\bar n}\frac{G^{(k)}_n(\xi)}{[n]_{q_k}^!}(x-\xi)_{q_k,n}-
    \sum_{n=0}^{\bar n}\frac{G_n(\xi)}{n!}(x-\xi)^n\r\|_\xi(R),\r.\\
&&\ds\l.\sup_{n\geq \bar n}\l\|\frac{G^{(k)}_n(\xi)}{[n]_{q_k}^!}(x-\xi)_{q_k,n}-
    \frac{G_n(\xi)}{n!}(x-\xi)^n\r\|_\xi(R)\r)<\veps\ .
\end{array}$$
\end{proof}

\makeatletter\newcount\uffa\uffa=\c@subsection\makeatother
\section{Weak Frobenius structure over a disk}
\label{sec:frobenius}
\makeatletter\c@subsection=\uffa\makeatother

{\bf Warning.}
\emph{In chapters \ref{sec:frobenius} and \ref{sec:transfer}
we will assume that $|1-q|<|\pi|=|\pi_q|$: this implies that for any integer $n$
we have
$$|1-q^n|=|\log q^n|=|n||\log q|=|n||1-q|\ .$$}
\par
\emph{Moreover throughout all of chapters \ref{sec:frobenius} and \ref{sec:transfer}
we will consider only $q$-difference systems over $D=D(0,1^-)$, therefore we will suppress
everywhere the index ``${}_{0,1}$'': so we will write $t$ for $t_{0,1}$,
$\|~\|$ for $\|~\|_{0,1}$, $\chi(A,q)$ for $\chi_{0,1}(A,q)$ and so on.}

\medskip
In this chapter we study the action of the Frobenius map $x\longmapsto x^{p^\ell}$,
$\ell\in\Z_{>0}$, on $q$-difference systems. Namely, following \cite{Ctransfer}, under convenient assumptions,
we will construct a matrix $H(x)$, depending on the $q$-difference system
$Y(qx)=A(x)Y(x)$, such that $A_{[H]}(x)$ is a function of
$x^{p^\ell}$.
This means that
$$
Y(qx)=A_{[H]}(x)Y(x)
$$
is actually a $q^{p^\ell}$-difference system in the variable $x^{p^\ell}$.
The interest of such a construction is that it changes the generic radius of convergence.
In fact,
$$
\chi(A_{[H]}(x^{p^\ell}),q^{p^\ell})=\chi(A(x),q)^{p^\ell}\ .
$$
\par
The original Christol theorem in \cite{Ctransfer} for $p$-adic differential equation,
as well as its generalization to differential modules over an annulus, due to Christol-Dwork \cite{CD},
is an irreplaceable tool: we think that the theorem above is destined to play
an analogous role for $q$-difference equations.

\subsection{Frobenius action on $q$-difference systems}

Until now we have worked with analytic and meromorphic functions:
in the next sections we will restrict our attention to analytic elements.
An introduction to the theory of analytic elements can be found
for instance in \cite[\S6, 4]{ROB}.
We briefly  recall the definition:

\begin{defn}
The ring $E_0$ of the
\emph{analytic elements on the disk $D(0,1^-)$ (defined over $K$)}
is the completion of the subring of $K(x)$ of all rational function
not having poles in $D(0,1^-)$ with respect to the norm $\|~\|$.
\par
We denote by $E_0^\p$ its quotient field.
\end{defn}

\begin{rmk}
\hfill\par\smallskip\noindent
1) Of course one can define analytic elements over any disk, but in this chapter we will
deal with disks centered at $0$. Therefore, by rescaling, we will always consider analytic
elements over $D(0,1^-)$.
\par\smallskip\noindent
2) In the sequel we will essentially use two properties of analytic elements, namely:
\begin{itemize}
\setlength{\itemsep}{-2pt}
\item By taking Taylor expansion at zero, one can identify analytic elements with bounded analytic functions over
$D(0,1^-)$: this ring embedding is an isometry with respect to $\|~\|$.
\item Any element of $E_0$ has a finite number of zeros in $D(0,1^-)$ (cf. \cite[IV, 5.2]{DGS}), hence any element
of $E_0^\p$ has a finite number of zeros and poles in $D(0,1^-)$.
\end{itemize}
\end{rmk}

The following theorem concerns a class ${\mathcal H}_q^\ell$, $\ell\in\Z_{>0}$,
of $q$-difference systems $Y(qx)=A(x)Y(x)$
satisfying the properties:
\begin{enumerate}
\setlength{\itemsep}{-2pt}
\item $A(x)\in Gl_\mu(E_0^\p)$;
\item $A(x)$ is analytic at $0$ and $A(0)\in Gl_\mu(K)$;
\item $Y(qx)=A(x)Y(x)$ has at worst apparent singularities in $D(0,1^-)\smallsetminus\{0\}$;
\item $\chi(A,q)>|\pi_q|^{\frac{1}{p^{\ell-1}}}$.
\end{enumerate}

We will say that a system in ${\mathcal H}_q^\ell$ is in \emph{normal form} if moreover
it satisfies the conditions 1 through 3 of the following proposition
(cf. \S\ref{subsec:step1} for the proof):

\begin{prop}\label{prop:normalform}
Let $Y(qx)=A(x)Y(x)$ be in ${\mathcal H}_q^\ell$.
Then there exists $U\in GL_\mu(K(x))$ such that the $q$-difference system associated to
$A_{[U]}(x)$ is in ${\mathcal H}_q^\ell$ and moreover:
\begin{enumerate}
\setlength{\itemsep}{-2pt}
\item $A_{[U]}(x)\in Gl_\mu(E_0)$;
\item any two eigenvalues $\la_1$, $\la_2$ of $A_{[U]}(0)$ satisfy  either $\la_1=\la_2$ or $\la_1\not\in \la_2 q^\Z$;
\item any eigenvalue $\la$ of $A_{[U]}(0)$ satisfies $\l|\frac{\la-1}{q-1}\r|<p^{1-\ell}$.
\end{enumerate}
\end{prop}

\noindent
Finally we are able to state the main theorem:

\begin{thm}\label{thm:frobenius}
Let us consider a $q$-difference system in ${\mathcal H}_q^\ell$
in its normal form
$$
Y(qx)=A(x)Y(x)\ .
$$
Then there exists a matrix $H(x)\in Gl_\mu(E^\p_0)\cap M_{\mu\times\mu}(E_0)$ such that
\begin{eqnarray}
&\hskip 20 pt\hbox to 11 true cm{$H(0)\in Gl_\mu(K)$ and $\l|H(0)\r|=\l|H(0)^{-1}\r|=1$;\hfill}\label{eq:uno}\\
&\hskip 20 pt\hbox to 11 true cm{$A_{[H]}(x)=F(x^{p^\ell})$;\hfill}\label{eq:due}\\
&\hskip 20 pt\hbox to 11 true cm{the eigenvalues of $A(0)$ and $F(0)$ coincide.\hfill}\label{eq:tre}
\end{eqnarray}
Moreover the $q^{p^\ell}$-difference system $V(q^{p^\ell}X)=F(X)V(X)$, with $X=x^{p^\ell}$,
has the following properties:
\begin{eqnarray}
&\hskip 20 pt\hbox to 11 true cm{$\chi(F,q^{p^\ell})=\chi(A,q)^{p^\ell}$;\hfill}\label{eq:quattro}\\
&\hskip 20 pt\hbox to 11 true cm{$F(X)\in Gl_\mu(E_0)$.\hfill}\label{eq:cinque}
\end{eqnarray}
\end{thm}

\begin{rmk}
It follows by (\ref{lemma:transferevidente}) that
a system $Y(qx)=A(x)Y(x)$ with coefficients in $E_0$
has only ordinary $q$-orbits in $D(0,1^-)\smallsetminus\{0\}$
if and only if there exists $n\geq 0$ such that
$x^nA(x)\in Gl_\mu(E_0)$. Therefore if $A(x)\in Gl_\mu(E_0)$ then
$Y(qx)=A(x)Y(x)$ has only ordinary $q$-orbits in $D(0,1^-)\smallsetminus\{0\}$ and $0$
is a regular singularity (cf. Appendix \ref{sec:regsing}).
Conversely, if $A(x)$ is analytic at $0$ with $A(0)\in Gl_\mu(K)$ and
$Y(qx)=A(x)Y(x)$ has only ordinary $q$-orbits in $D(0,1^-)\smallsetminus\{0\}$
then $A(x)\in Gl_\mu(E_0)$.
\end{rmk}

The proof of the above theorem, which follows the proof of Christol \cite{Ctransfer},
is the goal of \S\ref{subsec:prooffrobenius}.

\subsection{Normal form for a system in ${\mathcal H}_q^\ell$}
\label{subsec:step1}

Proposition \ref{prop:normalform} is a consequence of the following lemma:

\begin{lemma}\label{lemma:eigenvalue}
If $\chi(A,q)>|\pi_q|^{\frac{1}{p^{\ell-1}}}$, then for any eigenvalue $\la$
of $A(0)$ we have:
$$
{\rm dist}(\la, q^{\Z_p})=\inf_{\a\in\Z_p}|\la-q^\a|\leq p^{1-\ell}|q-1|\ ,
$$
where $q^\a$, for $\a\in\Z_p$, is defined as the sum of the binomial series
$$
q^\a=\sum_{n\geq 0}\frac{\a(\a-1)\cdots(\a-n+1)}{n!}(q-1)^n\ .
$$
\end{lemma}

\noindent
In fact:

\medskip\noindent
{\bf Proof of (\ref{prop:normalform}).}
There exists a \emph{shearing transformation} $U_1(x)$ with coefficients in $K\l[x,\frac{1}{x}\r]$,
(cf. Appendix \ref{sec:regsing}, in particular (\ref{prop:shearing}))
such that the eigenvalues of $A_{[U_1]}(0)$ are the eigenvalues of $A(0)$, multiplied by
chosen powers of $q$.
Hence by the previous lemma we can assume that any two eigenvalues $\la_1$, $\la_2$ of $A_{[U_1]}(0)$ satisfy the conditions:
\begin{trivlist}
\item 1) either $\la_1=\la_2$ or $\la_1\not\in \la_2 q^\Z$;
\item 2) $\ds\l|\frac{\la_1-1}{q-1}\r|<p^{1-\ell}$.
\end{trivlist}
Then by (\ref{prop:rem}) there exists a matrix $U_2(x)$, with coefficients in $K(x)$ and analytic at $0$, such that
$A_{[U_2U_1]}(x)$ has only ordinary orbits in $D\smallsetminus\{0\}$.
\par
Since $A_{[U_2U_1]}(x)$ is analytic at $0$ and $Y(qx)=A_{[U_2U_1]}(x)Y(x)$
has only ordinary orbits in $D\smallsetminus\{0\}$, we conclude that $A_{[U_2U_1]}(x)\in Gl_\mu(E_0)$.
Moreover
the eigenvalues of $A_{[U_1]}(0)$ and $A_{[U_2U_1]}(0)$ coincide, since $U_2(0)\in Gl_\mu(K)$. Therefore it is enough
to set $U(x)=U_2(x)U_1(x)$ to finish the proof.
\hfill$\square$

\medskip\noindent
Let us prove lemma \ref{lemma:eigenvalue}:

\medskip\noindent
{\bf Proof of (\ref{lemma:eigenvalue}).}
Write the system
$$
Y(qx)=A(x)Y(x)
\leqno{({\mathcal S})_q}
$$
in the form
\begin{equation}\label{eq:Gcorsivo}
\frac{\dq^n}{[n]_q^!}={\mathcal G}_n(x)Y(x)\ ,
\hbox{ for all $n\geq 0$.}
\end{equation}
Then
$$
{\mathcal G}_1(x)=\frac{1}{x}\frac{A(0)-{\mathbb I}_\mu}{q-1}+
\hbox{higher order terms.}
$$
By (\ref{eq:leibnizquattro}) and (\ref{eq:binomio}),
this formula generalizes to any $n\geq 1$ in the following way:
\begin{equation}\label{eq:poloG}
{\mathcal G}_n(x)
=\frac{q^{-\frac{n(n-1)}{2}}}{x^n}
    \frac{\l(A(0)-1\r)\l(A(0)-q\r)\cdots\l(A(0)-q^{n-1}\r)}
    {(q-1)(q^2-1)\cdots(q^n-1)}+\hbox{h.o.t.}
\end{equation}
We set
$$
\widetilde{\mathcal G}_n=
\frac{\l(A(0)-1\r)\l(A(0)-q\r)\cdots\l(A(0)-q^{n-1}\r)}
{(q-1)(q^2-1)\cdots(q^n-1)}\ .
$$
Since
$$
\l|\widetilde{\mathcal G}_n\r|\leq \l\|x^n{\mathcal G}_n(x)\r\|
$$
the series $g(x)=\sum_{n\geq 1}\widetilde{\mathcal G}_n x^n$
converges for $|x|<\chi(A,q)$, with $\chi(A,q)>|\pi_q|^{\frac{1}{p^{\ell-1}}}$.
\par
Let $C\in Gl_\mu(K)$ be a constant matrix such that $C^{-1}A(0)C$ is a
matrix in the Jordan normal form, then
$$
C^{-1}g(x)C=\sum_{n\geq 1}\frac{\l(C^{-1}A(0)C-1\r)\l(C^{-1}A(0)C-q\r)\cdots\l(C^{-1}A(0)C-q^{n-1}\r)}
{(q-1)(q^2-1)\cdots(q^n-1)}x^n\ .
$$
Hence we conclude that the eigenvalues of $A(0)$ satisfy the desired inequality
by applying the following result:

\begin{lemma}\label{lemma:eigenvalue2}
Let $\la\in K$ and let $\Lambda(x)=\sum_{n\geq 1}\la_nx^n$, with
$$
\la_n=\frac{(\la-1)\l(\la-q\r)\cdots\l(\la-q^{n-1}\r)}{(q-1)(q^2-1)\cdots(q^n-1)}\ .
$$
Then:
\par\noindent
1) The series $\Lambda(x)$ converges for $|x|<1$ if and only
if $\la\in q^{\Z_p}$, i.e. if and only if ${\rm dist}(\la,q^{\Z_p})=0$.
\par\noindent
2) If ${\rm dist}(\la,q^{\Z_p})=|p|^{k+\veps}|q-1|$,
with $k\in\Z_{\geq 0}$ and $\veps\in[0,1)$, then $\Lambda(x)$ converges for
$|x|<|p|^{-\frac{1}{p^k}\l(\frac{\veps}{p}-\frac{1}{p-1}\r)}$.
\par\noindent
3) If the series $\Lambda(x)$ converges for $|x|\leq |\pi|^{\frac{1}{p^{\ell-1}}}$, then
${\rm dist}(\la,q^{\Z_p})\leq p^{1-\ell}|q-1|$.
\end{lemma}

\begin{proof}\hfill
\par\smallskip\noindent
1) Let $\a\in\Z_p$ be such that $\la=q^\a$.
If $n_k$ is a sequence of integers such that $n_k\to\a$, then
the binomial series $x^{n_k}$ tends uniformly to $x^\a$ over any
closed disk $D(1,R^+)$, with $R<1$ (cf. \cite[IV, 5.4 and IV, \S7]{DGS}).
Hence, since $\l|q^n-1\r|\leq |q-1|<|\pi|$ for any $n\in\Z$, by the density of $\Z$ in $\Z_p$,
we have also $\l|q^{\a-n}-1\r|<|\pi|$.
We deduce  from the $p$-adic properties of logarithm that
\begin{equation}\label{eq:log}
\begin{array}{rcl}|\la_n|
&=&\ds\l|\frac{(q^\a-1)\l(q^\a q^{-1}-1\r)\cdots\l(q^\a q^{1-n}-1\r)}
    {(q-1)(q^2-1)\cdots(q^n-1)}\r|\\ \\
&=&\ds\frac{|\log q^\a||\log(q^\a q^{-1})|\cdots|\log(q^\a q^{1-n})|}
    {|\log q|^n|n!|}\\ \\
&=&\ds\frac{1}{|n!|}\l|\frac{\log q^\a}{\log q}\r|\l|\frac{\log q^\a}{\log q}-1\r|
    \cdots\l|\frac{\log q^\a}{\log q}-(n-1)\r|\ .
\end{array}
\end{equation}
Since $\log q^n=n\log q$ for any $n$ and $\Z$ is dense in $\Z_p$,
one verifies that $\frac{\log q^\a}{\log q}=\a$. We conclude that
the radius of convergence of $\Lambda(x)$ is equal to the radius of convergence of
the binomial series $\sum_{n\geq 0}{\a\choose n}x^n$, hence is equal to $1$ \cite[IV, 7.5]{DGS}.
\par
Let us prove the opposite implication.
Remark that
$$
\la_n=\frac{1}{[n]_q^!}
\frac{\la-1}{q-1}\l(\frac{\la-1}{q-1}-[1]_q\r)\cdots\l(\frac{\la-1}{q-1}-[n-1]_q\r)\ .
$$
Hence if $\l|\frac{\la-1}{q-1}\r|>1$ we have
$$
\l|\la_n\r|=\l|\frac{\la-1}{q-1}\r|^n\frac{1}{|n!|}
$$
and the radius of convergence of $\Lambda(x)$ would be $\l|\pi\frac{q-1}{\la-1}\r|<1$.
We conclude that $\Lambda(x)$ convergent for $|x|<1$ implies $|\la-1|\leq|q-1|<|\pi|$.
Suppose that $\Lambda(x)$ is convergent for $|x|<1$.
As in (\ref{eq:log}) we obtain:
$$
|\la_n|=\frac{1}{|n!|}\l|\frac{\log \la}{\log q}\r|\l|\frac{\log \la}{\log q}-1\r|
    \cdots\l|\frac{\log \la}{\log q}-(n-1)\r|\ .
$$
We conclude that $\a=\frac{\log\la}{\log q}\in\Z_p$ since
the series $\sum_{n\geq 0}{\log\la/\log q\choose n}x^n$ converges for
$|x|<1$ \cite[IV, 7.5]{DGS}.
Then $\log\la=\log q^\a$, with $|\la-1|,|q^\a-1|<|\pi|$, and therefore
$$
\la=\exp\log\la=\exp\log q^\a=q^\a\ .
$$
\par\smallskip\noindent
2)
Since $\Z_p$ is compact there exists $\a\in\Z_p$ such that $|\la-q^\a|=|q-1||p|^{k+\veps}$.
Let $s$ be a nonnegative integer such that
$|q^s-q^\a|<|q-1||p|^{k+\veps}$; then $|\la-q^s|=|q-1||p|^{k+\veps}={\rm dist}(\la,q^{\Z_p})$.
Hence we can choose $\a=s\in\Z_{\geq0}$.
\par
Let us suppose $s=0$.
Then for any $m\in\Z$ we have:
$$
\l|\frac{\la-q^m}{q-1}\r|
=\l|\frac{\la-1}{q-1}-\frac{q^m-1}{q-1}\r|
\geq\l|\frac{\la-1}{q-1}\r|
$$
We obtain the estimate (cf. for instance \cite[2.2]{BCII} or \cite[IV, 7.3]{DGS}
for the analogous classical estimate):
$$
{\rm ord}_p ([n]_q^!\la_n)=\sum_{i=1}^k\l(1+\l[\frac{n-1}{p^i}\r]\r)+
\veps\l[\frac{n-1}{p^{k+1}}\r]
={\rm ord}_p {\frac{\la-1}{q-1}\choose n}\ ,
$$
with
$$
{\frac{\la-1}{q-1}\choose n}=\frac{1}{n!}\frac{\la-1}{q-1}\l(\frac{\la-1}{q-1}-1\r)\cdots\l(\frac{\la-1}{q-1}-n+1\r)\ .
$$
Since $|[n]_q^!|=|n!|$, the radius of convergence of $\Lambda(x)$ is the same as the radius of convergence of the
binomial series $\sum_{n\geq 0}{\frac{\la-1}{q-1}\choose n}x^n$,
namely
$p^{\frac{1}{p^k}\l(\frac{\veps}{p}-\frac{1}{p-1}\r)}$
(see \cite[2.2]{BCII} or \cite[IV, 7.3]{DGS} for the explicit calculation).
\par\smallskip\noindent
3) If $\la\in q^{\Z_p}$ the assertion follows immediately from 1), hence it is enough to consider
the case $\la\not\in q^{\Z_p}$.
Moreover, if $\l|\frac{\la-1}{q-1}\r|>1$ the series $\Lambda(x)$ has radius of convergence
$\l|\pi\frac{q-1}{\la-1}\r|<|\pi|^{\frac{1}{p^{\ell-1}}}$, therefore
$\la$ necessarily satisfies ${\rm dist}(\la,q^{\Z_p})\leq |q-1|$.
Let
$$
\frac{{\rm dist}(\la,q^{\Z_p})}{|q-1|}=\l|\frac{\la-1}{q-1}\r|=|p|^{k+\veps}\ ,
$$
with $k\in\Z_{\geq 0}$ and $\veps\in[0,1)$.
We finish the proof by observing that
$|\pi|^{\frac{1}{p^{\ell-1}}}\leq p^{\frac{1}{p^k}\l(\frac{\veps}{p}-\frac{1}{p-1}\r)}$
implies that $\ell-1<k+\veps$ (cf. again \cite[2.2]{BCII} or \cite[IV, 7.3]{DGS}).
\par
Finally if $s>0$, then for any $n>s$ we have:
$$
\la_n=\la_s
\frac{(q^{-s}\la-1)\l(q^{-s}\la-q\r)\cdots\l(q^{-s}\la-q^{n-s-1}\r)}
{(q-1)(q^2-1)\cdots(q^{n-s-1}-1)}\ .
$$
By taking $q^{-s}\la$ instead of $\la$, we reduce to the case $s=0$.
\end{proof}

\subsection{Proof of (\ref{thm:frobenius})}
\label{subsec:prooffrobenius}

The proof of (\ref{thm:frobenius}) is divided into two steps:
\begin{trivlist}
\item {\it Step 1}: construction of the matrix $H(x)$ and the proof of (\ref{eq:uno}).
\item {\it Step 2}: proof of (\ref{eq:due}) through (\ref{eq:cinque}).
\end{trivlist}
They are handled respectively in (\ref{parag:step2}) and
(\ref{parag:step3}).

\begin{parag}
{\bf Construction of the matrix $H(x)$ and proof of (\ref{eq:uno}).}
\label{parag:step2}

\par\medskip\noindent
By assumption (cf. (\ref{eq:Gcorsivo}) for
notation) the solution
$
Y(\xi,x)=\sum_{n\geq 0}{\mathcal G}_n(\xi)(x-\xi)_{q,n}
$
of $Y(qx)=A(x)Y(x)$ is defined for any
$\xi\in (D(0,1^-)\smallsetminus\{0\})\cup D(t,1^-)$ and
converges for $|x-\xi|<\chi(A,q)$.
We set:
\begin{equation}
\begin{array}{rcl}
H(x)
&:=&\ds\frac{1}{p^\ell}\sum_{\zeta^{p^\ell}=1}Y(x,\zeta x)
=\ds\frac{1}{p^\ell}\sum_{n\geq 0}{\mathcal G}_n(x)x^n\sum_{\zeta^{p^\ell}=1}(\zeta-1)_{q,n}\ .
\end{array}
\end{equation}
Since ${\mathcal G}_n(x)x^n\in M_{\mu\times\mu}({E_0})$, with
$$
\liminf_{n\rightarrow\infty}\l\|{\mathcal G}_n(x)x^n\r\|^{-1/n}=\chi(A,q),
$$
and
$$
\begin{array}{rcl}
|(\zeta-1)_{q,n}|
&=&\ds\prod_{i=0}^n|\zeta-1+1-q^i|\\
&\leq&\ds\prod_{i=0}^n\sup\l(|\zeta-1|, |i||1-q|\r)\\
&\leq&\ds|\pi|^{\frac{n}{p^{\ell-1}}}<\chi(A,q)^n\ ,
\end{array}
$$
the matrix $H(x)$ is a well-defined element of $M_{\mu\times\mu}({E_0})$
and it makes sense to evaluate $H(x)$ at $0$.
\par
Notice that $x^n{\mathcal G}_n(x)$ has no pole at zero. In the notation of the
proof of lemma \ref{lemma:eigenvalue} its constant term is equal to $q^{-\frac{n(n-1)}{2}}\widetilde{\mathcal G}_n$.
Hence we have
$$
H(0)=\frac{1}{p^\ell}\sum_{n\geq 0}q^{-\frac{n(n-1)}{2}}
\widetilde{\mathcal G}_n\sum_{\zeta^{p^\ell}=1}(\zeta-1)_{q,n}
$$
By (\ref{lemma:eigenvalue2}), $H(0)$ converges in $K$.
Since (cf. \cite[IV, 7.3]{DGS})
$$
\frac{1}{p^\ell}\sum_{n\geq 0}\sum_{\zeta^{p^\ell}=1}
{\widetilde{\mathcal G}_1\choose n}(\zeta-1)^n={\mathbb I}_\mu
$$
and
$$
\l|H(0)-{\mathbb I}_\mu\r|
=\l|\frac{1}{p^\ell}\sum_{n\geq 0}\sum_{\zeta^{p^\ell}=1}
      \l(q^{-\frac{n(n-1)}{2}}\widetilde{\mathcal G}_n(\zeta-1)_{q,n}
      -{\widetilde{\mathcal G}_1\choose n}(\zeta-1)^n\r)\r|<1\ ,
$$
we conclude that $|H(0)|=1$. Since moreover $K$ is complete, this also proves that $H(0)$
has an inverse, which has norm $1$.
\end{parag}

\begin{parag}
{\bf The matrix $H(x)$ satisfies (\ref{eq:due}) through (\ref{eq:cinque}).}
\label{parag:step3}

\par\medskip\noindent
By assumption the solution $Y(t,x)$ converges for $|x-t|<\chi(A,q)$.
It follows from \cite[VI, \S6]{DGS} that the matrix
$$
V(x)=\frac{1}{p^\ell}\sum_{z^{p^\ell}=x}Y(t,z)
$$
is analytic over the disk $D(t^{p^\ell},\chi(A,q)^{p^{\ell}})$.
The formula (\ref{eq:solformali}) implies that
$$
\begin{array}{rcl}
V(x^{p^\ell})
&=&\ds\frac{1}{p^\ell}\sum_{z^{p^\ell}=x^{p^\ell}}Y(t,z)\\
&=&\ds\frac{1}{p^\ell}\sum_{\zeta^{p^\ell}=1}Y(t,\zeta x)
    \frac{1}{p^\ell}\sum_{\zeta^{p^\ell}=1}Y(x,\zeta x)Y(t,x)\\
&=&\ds H(x)Y(t,x)\ .
\end{array}
$$
Hence the matrix
$$
F(x)=V(q^{p^\ell}x)V(x)^{-1}
$$
satisfies (\ref{eq:due}). In fact,
$$
\begin{array}{rcl}
F(x^{p^\ell})
&=&\ds V(q^{p^\ell}x^{p^\ell})V(x^{p^\ell})^{-1}\\
&=&\ds H(qx)Y(t,qx)Y(t,x)^{-1}H(x)^{-1}\\
&=&\ds H(qx)A(x)H(x)^{-1}\ .
\end{array}
$$
Then (\ref{eq:tre}) follows immediately from (\ref{eq:uno}).

\medskip
Let $X=x^{p^\ell}$. Then $V(X)$ is solution to
$$
V(q^{p^\ell}X)=F(X)V(X)\ .
$$
Since $V(X)$ is analytic around $t^{p^\ell}$ and
$V(t^{p^\ell})=H(t)$, the generic radius $\chi(F,q^{p^\ell})$ at $t^{p^\ell}$
coincides with the radius of convergence of $H(t)^{-1}V(X)$. Hence
$$
\chi(F,q^{p^\ell})=\chi(A,q)^{p^\ell}\ .
$$
Finally let $\xi\in D(0,1^-)$, $\xi\neq 0$. Similar calculations to the ones we have
worked out for $V(x)$ shows that an analytic fundamental solution of
$V(q^{p^\ell}X)=F(X)V(X)$ at $\xi^{p^\ell}$ is given by
$$
\frac{1}{p^\ell}\sum_{z^{p^\ell}=x^{p^\ell}}Y(\xi,z)\ .
$$
Since
$$
\begin{array}{rcl}
D(0,1^-)&\longrightarrow &D(0,1^-)\\
\xi     &\longmapsto     &\xi^{p^\ell}
\end{array}
$$
is a bijective map, this proves (\ref{eq:cinque}).
\end{parag}

\makeatletter\newcount\uffa\uffa=\c@subsection\makeatother
\section{Transfer theorems in regular singular disks}
\label{sec:transfer}
\makeatletter\c@subsection=\uffa\makeatother

In \S\ref{subsec:consequences} we have seen how the generic radius of convergence
is linked to the radius of convergence of analytic solutions at ordinary
$q$-orbits.
By (\ref{prop:rem}) we are able to deduce estimates also for the radius of meromorphy
of solutions at
apparent and trivial singularities. In this chapter we are going to consider
a $q$-difference system having a regular singularity at $0$
(cf. Appendix \ref{sec:regsing} for a summary of basic properties of regular singular
$q$-difference equations). In particular we are going to prove a
$q$-difference version of the Christol-Andr\'e-Baldassarri-Chiarellotto theorem
(cf. \cite{Ctransfer}, \cite{Atransfer}, \cite{BCI} and \cite{BCII}).

\subsection{An analogue of Christol's theorem}

A $q$-difference system
$$
Y(qx)=A(x)Y(x)
$$
with coefficients in $E_0^\p$
is \emph{regular singular at $0$} if it is regular singular regarded as a
$q$-difference system with coefficients
in the field $K((x))$ via the canonical immersion $E^\p_0\hookrightarrow K((x))$,
i.e. if there exists $U(x)\in Gl_\mu(K((x)))$ such that $A_{[U]}(x)\in Gl_\mu(K\[[x\]])$.
\par
Let $r(U(x))$ be the radius of convergence of the matrix $x^N U(x)$, for $N$ large enough to have
$x^NU(x)\in Gl_\mu(K\[[x\]])$.
Very similarly to the differential case,
the theorem below establishes an estimate of $r(U(x))$ with respect to the generic radius of convergence and
a certain number attached to the eigenvalues of $A_{[U]}(0)$:

\begin{defn}
We call the \emph{$q$-type of $\a\in K$}, and we write ${\rm type}_q(\a)$,
the radius of convergence of
$$
\sum_{n\geq 0,\a\neq q^n}\frac{1-q}{1-q^n\a}x^n\ .
$$
\end{defn}

\begin{rmk}\label{rmk:type}
In the $p$-adic theory of differential equations one calls the type of a number $\a\in K$ the radius
of convergence ${\rm type}(\a)$ of the series
$$
\sum_{n\geq 0,n\neq\a}\frac{x^n}{n-\a}\ .
$$
The definition above is an analogue of this notion
(see Appendix \ref{sec:type}, \S\ref{subsec:type}, for an estimate of the $q$-type
in terms of the classical type and some general properties of the $q$-type).
\end{rmk}

\medskip\noindent
The following is a $q$-analogue of Christol theorem \cite{Ctransfer}. It is a
transfer theorem for $q$-difference systems $Y(qx)=A(x)Y(x)$
such that $\chi(A,q)=1$. Statements concerning the general situation can be found
in \S\ref{subsec:generalstatements}.

\begin{thm}\label{thm:transfer}
Let us consider a system $Y(qx)=A(x)Y(x)$ such that
\begin{trivlist}
\item 1) $A(x)\in Gl_\mu(E_0^\p)$;
\item 2) the system $Y(qx)=A(x)Y(x)$ has only apparent singularities
in $D(0,1^-)^\smallsetminus\{0\}$;
\item 3) the system $Y(qx)=A(x)Y(x)$ has a regular singulairty at $0$,
i.e. there exists $U(x)\in Gl_\mu(K((x)))$
such that $A_{[U]}\in M_{\mu\times\mu}(K)$;
\item 4) $\chi(A,q)=1$.
\end{trivlist}
Then
\begin{equation}
r(U(x))\geq\prod{\rm type}_q\l(\a\be^{-1}\r)\ ,
\end{equation}
where the product is taken over all the (ordered) couples of eigenvalues
$\a,\be$ of $A_{[U]}$.
\end{thm}

The proof of this theorem, as for the Christol-Andr\'e-Baldassarri-Chiarellotto theorem,
which inspired it, relies on the existence of the weak Frobenius structure and
it is quite long. Sections \S\ref{subsec:rough} and \S\ref{subsec:sharp} are devoted to the
proof of some preliminary estimates of $r(U(x))$.
As a corollary, in \S\ref{subsec:generalstatements} we prove the transfer theorem
\ref{thm:transfer} plus some more general statements, without any assumption on $\chi(A,q)$.

\subsection{A first rough estimate}
\label{subsec:rough}

In this section we prove a first estimate, which is not very sharp, but
crucial for the proof of (\ref{thm:transfer}):

\begin{prop}\label{prop:firstestimate}
Let $Y(qx)=A(x)Y(x)$ be a system such that
\begin{trivlist}
\item 1) $A(x)\in Gl_\mu(E_0)$;
\item 2) the system $Y(qx)=A(x)Y(x)$ has only ordinary $q$-orbits
in $D(0,1^-)^\smallsetminus\{0\}$;
\item 3) any two eigenvalues $\a,\be$ of $A(0)$ satisfy either $\a=\be$ or
$\a\be^{-1}\not\in q^\Z$.
\end{trivlist}
Then the matrix $U(x)\in Gl_\mu(K\[[x\]])$ satisfying
$A_{[U]}(x)=A(0)$ is such that
$$
r(U(x))\geq
\inf\l(|\pi|,\chi(A,q)\r)^{\mu^2}\l(\det A(0)\r)^\mu\prod_{\a,\be\hbox{ e.v. of }A(0)}
{\rm type}_q\l(\a\be^{-1}\r)\ .
$$
\end{prop}

\begin{rmk}
Observe that under the assumption of the proposition above, the matrix $U(x)$ always exists
(cf. Appendix \ref{sec:regsing}, (\ref{lemma:uniforme})).
\end{rmk}

\begin{proof}
The proof is divided into two steps:

\par\smallskip\noindent
{\it Step 1. One can assume that $Y(qx)=A(x)Y(x)$ satisfies
\begin{equation}\label{eq:hyp4}
\l\|\frac{A(x)-{\mathbb I}_\mu}{(q-1)x}\r\|\leq \sup\l(1,\frac{|\pi|}{\chi(A,q)}\r)\ .
\end{equation}}

\smallskip\noindent
Let us set $\chi=\chi(A,q)$ to simplify notation. By the cyclic vector lemma
there exists a matrix $H_1(x)\in Gl_\mu(E_0^\p)$ such that the matrix $A_{[H_1]}(x)$
is in the form
$$
A_{[H_1]}(x)=\l(\begin{array}{c|c}
0& \\
\vdots &{\mathbb I}_{\mu-1}\\
0& \\
\hline\\
a_0(x)&a_1(x)\dots a_{\mu-1}(x)
\end{array}\r)\ .
$$
An analogue of Fuchs' theory for $q$-difference equation
(cf. Appendix \ref{sec:regsing}, \S\ref{subsec:qequation}) assures that
$A_{[H_1]}(x)$ is analytic at $0$ and $A_{[H_1]}(0)$ is an
invertible constant matrix.
Consider the matrix $\widetilde H$:
$$
\widetilde H=\l(a_{i,j}\r)_{i,j=0,\dots,\mu-1}\ ,
\hbox{ with }
a_{i,j}=\begin{cases}
\ds\frac{(-1)^j}{(q-1)^i}{i\choose j}& \hbox{ if }j\leq i;\\ \\
0&\hbox{ otherwise}.
\end{cases}
$$
Then we have (cf. (\ref{eq:leibnizcinque}))
$$
xG_{[\tilde HH_1]}=\frac{A_{[\tilde HH_1]}(x)-{\mathbb I}_\mu}{q-1}
\l(\begin{array}{c|c}
0& \\
\vdots &{\mathbb I}_{\mu-1}\\
0& \\
\hline\\
b_0(x)&b_1(x)\dots b_{\mu-1}(x)
\end{array}\r)\ .
$$
Remark that $xG_{[\tilde HH_1]}$ is analytic at $0$ and $A_{[\tilde HH_1]}(0)\in Gl_\mu(K)$, since
$A_{[H_1]}(x)$ has the same properties.
\par
By proposition
\ref{prop:dworkfrobenius}, if $\chi\geq|\pi|$, then
$$
\l\|G_{[\tilde HH_1]}\r\|\leq 1\ .
$$
On the other hand, if $\chi<|\pi|$, we immediately deduce from
(\ref{lemma:genradius}) and (\ref{prop:dworkfrobenius}) that
$$
\chi=\frac{|\pi|}{\sup_i\l\|b_i(x)\r\|^{1/(\mu-i)}}\ .
$$
Since $K$ is algebraically closed, we can find $\ga\in K$ such that
$$
|\ga|=\sup_i\l\|b_i(x)\r\|^{1/(\mu-i)}\ ;
$$
then an explicit calculation shows that the
diagonal gauge transformation matrix
$$
H_2=\begin{pmatrix}1&     &      & \\
                   &1/\ga&      & \\
               &     &\ddots& \\
               &     &      &1/\ga^{\mu-1}
    \end{pmatrix}
$$
satisfies
$$
\l\|G_{[H_2\tilde HH_1]}(x)\r\|\leq \frac{|\pi|}{\chi}\ .
$$
We set $H(x)=\widetilde HH_1(x)$ if $\chi\geq|\pi|$ and $H(x)=H_2\widetilde HH_1(x)$ otherwise.
In both cases $H(x)\in Gl_\mu(E_0^\p)$,
$A_{[H]}(x)$ is analytic at $0$  with $A_{[H]}(0)\in Gl_\mu(K)$, since
$A_{[\tilde HH_1]}(x)$ has the same properties, and
$$
\l\|G_{[H]}(x)\r\|\leq \sup\l(1,\frac{|\pi|}{\chi}\r)\ .
$$
We proceed as in (\ref{prop:normalform}):
applying successively a unimodular shearing transformation $H_4(x)\in
Gl_\mu\l(K\l[x,\frac{1}{x}\r]\r)$,
constructed as in (\ref{prop:shearing}), and a unimodular gauge transformation $H_3(x)\in Gl_\mu(K(x))$,
constructed as in
(\ref{prop:rem}) to remove apparent singularities, we obtain a matrix
$A_{[H_4H_3H]}(x)$ satisfying hypothesis 1), 2) and 3) plus the inequality:
$$
\l\|\frac{A_{[H_4H_3H]}(x)-{\mathbb I}_\mu}{(q-1)x}\r\|\leq\sup\l(1,\frac{|\pi|}{\chi(A,q)}\r)\
$$
Set $Q(x)=H_4(x)H_3(x)H(x)$ and $B(x)=A_{[Q]}(x)$.
One observes that
\begin{trivlist}
\item - since both $A(x)$ and $B(x)$ satisfy hypothesis 1) and 2), necessarily
$Q(x)$ and $Q(x)^{-1}$ have no poles in $D(0,1^-)\smallsetminus\{0\}$;
\item - $B_{[UQ^{-1}]}(x)=A(0)$.
\end{trivlist}
Once again there exists a shearing transformation
$P(x)\in Gl_\mu\l(K\l[x,\frac{1}{x}\r]\r)$ such that $B(0)=B_{[PUQ^{-1}]}(x)$.
Moreover $PUQ^{-1}$ has its coefficients in $K((x))$.
Set $V=PUQ^{-1}$.
Then the $q$-difference system $Y(qx)=B(x)Y(x)$ satisfies conditions
1) through 3), with $B_{[V]}(x)=B(0)$, plus the inequality (\ref{eq:hyp4}).
Since $r(V)=r(U)$, we conclude that it is enough to prove the
theorem assuming (\ref{eq:hyp4}).

\par\smallskip\noindent
{\it Step 2. Proof of ({\rm \ref{thm:transfer}}) assuming ({\rm \ref{eq:hyp4}}).}

\par\smallskip\noindent
Let us set
$$
U(x)=U_0+U_1x+U_2x^2+\dots
$$
and
$$
A(x)=A_0+A_1x+A_2x^2+\dots\ .
$$
Of course one can suppose $U_0={\mathbb I}_\mu$.
By assumption we have $U(qx)A(x)=A_0U(x)$.
Hence a direct calculation shows that
$$
q^mU_mA_0-A_0U_m=A_1U_{m-1}+A_2U_{m-2}+\dots+A_m\ .
$$
Consider the $K$-linear map
\begin{equation}
\begin{array}{rcl}
\Phi_{q^m,A_0}: M_{\mu\times\mu}(K) & \longrightarrow& M_{\mu\times\mu}(K)\\ \\
            M                   & \longmapsto    & q^m MA_0-A_0M
\end{array}\ .
\end{equation}
Observe that the eigenvalues of $\Phi_{q^m,A_0}$ are precisely of the form
$q^m\a-\be$, where $\a,\be$ are any two eigenvalues of $A_0$.
Therefore the operator $\Phi_{q^m,A_0}$ is invertible by hypothesis 3), which means that
$$
U_m=\Phi_{q^m,A_0}^{-1}\l(A_1U_{m-1}+A_2U_{m-2}+\dots+A_m\r)\ .
$$
We need to calculate the norm of $\Phi_{q^m,A_0}^{-1}$ as a $K$-linear operator.
Since $\Phi_{q^m,A_0}(M)=q^m M(A_0-{\mathbb I}_\mu)-M(q^m{\mathbb
I}_\mu-A_0)$ and
$$
\|A(x)-{\mathbb I}_\mu\|\leq|q-1|\sup\l(\l\|\frac{A(x)-{\mathbb I}_\mu}{(q-1)x}\r\|, 1\r)
\leq|q-1|\sup\l(\frac{|\pi|}{\chi},1\r)\ ,
$$
we conclude that
$$
\l\|\Phi_{q^m,A_0}^{-1}\r\|=
\frac{\l\|{\rm adj}\Phi_{q^m,A_0}\r\|}{\l\|\det \Phi_{q^m,A_0}\r\|}
\leq\frac{|q-1|^{\mu^2-1}\sup\l(\frac{|\pi|}{\chi}, 1\r)^{\mu^2-1}}
{\prod_{\a,\be\hbox{ e.v. of }A_0}|q^m\a-\be|}\ .
$$
Recursively we obtain the estimate
$$
\begin{array}{rcl}|U_m|
&\leq&\ds\l\|\Phi_{q^m,A_0}^{-1}\r\|\l(\sup_{i\geq 1}|A_i|\r)\l(\sup_{i=1,\dots, m-1}|U_i|\r)\\ \\
&\leq&\ds\prod_{i=1}^m\l\|\Phi_{q^m,A_0}^{-1}\r\|\|A(x)-{\mathbb I}_\mu\|^m\\ \\
&\leq&\ds\frac{\sup\l(\frac{|\pi|}{\chi}, 1\r)^{\mu^2m}}
    {|\det A_0|^{\mu m}}\l(\prod_{\a,\be\hbox{ e.v. of }A_0}
    \prod_{i=1}^m\l|\frac{1-q^m\a\be^{-1}}{1-q}\r|\r)^{-1}\ .
\end{array}
$$
We deduce from (\ref{cor:type}) that
$$
r(U(x))=\liminf_{m\rightarrow\infty}\l|U_m\r|^{-\frac{1}{m}}
\geq\inf\l(\l|\pi\r|,\chi\r)^{\mu^2}
|\det A_0|^\mu\prod_{\a,\be}{\rm type}_q(\a\be^{-1})\ .
$$

\end{proof}

For further reference, we point out that the Step 1 above
is actually a proof of the following statement:

\begin{lemma}\label{lemma:reduction}
Let $Y(qx)=A(x)Y(x)$ be a $q$-difference equation, with
$A(x)\in Gl_\mu(E^\p_0)$, having only apparent singularities in $D(0,1^{-})\smallsetminus\{0\}$
and a regular singularity at $0$.
We suppose that there exists $U(x)\in Gl_\mu(K((x)))$ satisfying
$A_{[U]}(x)\in Gl_\mu(K)$. Then there exists a gauge transformation matrix $H(x)\in Gl_\mu(E^\p_0)$
such that:
\begin{trivlist}
\item 1) $A_{[H]}(x)\in Gl_\mu(E_0)$;
\item 2) any two eigenvalue $\a,\be$ of $A_{[H]}(0)$ are either equal or $\a\be^{-1}\not\in q^\Z$;
\item 3) $r(U(x))=r(H(x)U(x))$\ .
\end{trivlist}
\end{lemma}

\subsection{A sharper estimate}
\label{subsec:sharp}

In this section we are going to the deduce from (\ref{prop:firstestimate})
a sharper estimate for $q$-difference system in ${\mathcal H}_q^\ell$
(cf. (\ref{prop:normalform})),
relying on the existence of the weak Frobenius structure.

\begin{prop}\label{prop:sharp}
Let $Y(qx)=A(x)Y(x)$ be a $q$-difference system in ${\mathcal H}_q^\ell$ in
its normal form and let $\chi(A,q)>|\pi|^{\frac{1}{p^{\ell-1}}}$.
The matrix $U(x)\in Gl_\mu(K\[[x\]])$ such that
$A_{[U]}=A(0)$ satisfies the inequality
$$
r(U(x))\geq
\inf\l(|\pi|^{\frac{1}{p^\ell}},\chi(A,q)\r)^{\mu^2}
\l(\det A(0)\r)^{\frac{\mu}{p^\ell}}\prod_{\a,\be\hbox{ e.v. of }A(0)}
{\rm type}_q\l(\a\be^{-1}\r)\ .
$$
\end{prop}

\begin{proof}
By theorem \ref{thm:frobenius} stating the existence of the weak Frobenius structure,
there exists a matrix $H(x)\in Gl_\mu(E^\p_0)\cap M_{\mu\times\mu}(E_0)$ such that:
\begin{trivlist}
\item 1) $H(0)\in Gl_\mu(K)$ and $\l|H(0)\r|=\l|H(0)^{-1}\r|=1$;
\item 2) $A_{[H]}(x)=F(x^{p^\ell})$;
\item 3) the eigenvalues of $A(0)$ and $F(0)$ coincide;
\item 4) the $q^{p^\ell}$-system $V(q^{p^\ell}X)=F(X)V(X)$, with $X=x^{p^\ell}$, has only ordinary $q$-orbits
    in $D(0,1^-)\smallsetminus\{0\}$;
\item 5) $\chi(F,q^{p^\ell})= \chi(A,q)^{p^\ell}$.
\end{trivlist}
By (\ref{lemma:uniforme}) there exists $V(X)\in Gl_\mu(K\[[X\]])$ such that
$F_{[V]}(X)=F(0)$. On the other hand we have also $F_{[UH^{-1}]}(x)=A(0)$.
It follows by (\ref{cor:uniforme}) that there exists a shearing transformation
$Q(x)\in Gl_\mu\l(K\l[x,\frac{1}{x}\r]\r)$ such that
$Q(x)U(x)=V(x)H(x)$. Hence $r(U(x))=r(V(x))=r(V(X))^{\frac{1}{p^\ell}}$.
\par
By (\ref{prop:firstestimate}) we have
$$
r(V(X))\geq \inf\l(|\pi|,\chi(A,q)^{p^\ell}\r)^{\mu^2}
\l(\det F(0)\r)^\mu\prod_{\a,\be\hbox{ e.v. of }F(0)}
{\rm type}_{q^{p^\ell}}\l(\a\be^{-1}\r)
$$
and hence
$$
r(U(x))\geq \inf\l(|\pi|^{\frac{1}{p^\ell}},\chi(A,q)\r)^{\mu^2}
\l(\det A(0)\r)^{\frac{\mu}{p^\ell}}\prod_{\a,\be\hbox{ e.v. of }A(0)}
{\rm type}_{q^{p^\ell}}\l(\a\be^{-1}\r)^{\frac{1}{p^\ell}}\ .
$$
We finish the proof by observing that
$$
{\rm type}_q(\a)=
\liminf_{n\rightarrow\infty}|1-q^n\a|^{1/n}
\leq\liminf_{n\rightarrow\infty}|1-q^{p^\ell n}\a|^{1/p^\ell n}
=\l({\rm type}_{q^{p^\ell}}\l(\a\r)\r)^{\frac{1}{p^\ell}}\ .
$$
\end{proof}

\subsection{More general statements}
\label{subsec:generalstatements}

Let $Y(qx)=A(x)Y(x)$ be a $q$-difference system such that
$A(x)\in Gl_\mu(E_0^\p)$, having only apparent singularities in
$D(0,1^-)\smallsetminus\{0\}$. We suppose that there exists
$U(x)\in Gl_\mu(K((x)))$ satisfying $A_{[U]}(x)\in Gl_\mu(K)$.
Then we have:

\begin{cor}\hfill
\par\smallskip\noindent
1) If $\chi(A,q)\leq|\pi|$.
Then
$$
r(U(x))\geq\chi(A,q)^{\mu^2}|\det A(0)|^{\mu}
\prod_{\a,\be\hbox{ e.v. of }A(0)}
{\rm type}_q\l(\a\be^{-1}\r)\ .
$$
\par\smallskip\noindent
2) Suppose $|\pi|<\chi(A,q)<1$ and let $\ell$ be a positive integer such that
$|\pi|^{\frac{1}{p^\ell}}\geq \chi(A,q)>|\pi|^{\frac{1}{p^{\ell-1}}}$.
Then
$$
r(U(x))\geq\chi(A,q)^{\mu^2}|\det A(0)|^{\frac{\mu}{p^\ell}}
\prod_{\a,\be\hbox{ e.v. of }A(0)}
{\rm type}_q\l(\a\be^{-1}\r)
$$
\par\smallskip\noindent
3) If $\chi(A,q)=1$.
Then
$$
r(U(x))\geq\prod_{\a,\be\hbox{ e.v. of }A(0)}
{\rm type}_q\l(\a\be^{-1}\r)
$$
\end{cor}

\begin{proof}
Statement 1) follows immediately from (\ref{lemma:reduction}) and (\ref{prop:firstestimate}),
while statement 2)
follows from (\ref{lemma:reduction}) and (\ref{prop:sharp}).
Let us prove 3).
By (\ref{lemma:reduction}) and (\ref{prop:sharp}), for any $\ell\in\Z_{\geq 0}$,
we have $\chi(A,q)>|\pi|^{\frac{1}{p^{\ell-1}}}$, hence
$$
r(U(x))\geq\l|\pi(1-q)\r|^{\frac{1}{p^\ell}}|\det A(0)|^{\frac{\mu}{p^\ell}}
\prod_{\a,\be\hbox{ e.v. of }A(0)}
{\rm type}_q\l(\a\be^{-1}\r)\ ,
\hbox{ for any $\ell>0$.}
$$
We finish by letting $\ell\rightarrow+\infty$.
\end{proof}


\makeatletter\newcount\uffa\uffa=\c@subsection\makeatother
\begin{appendix}
\section*{Appendix}

\section{Twisted Taylor expansion of $p$-adic analytic functions}
\label{sec:taylor}
\makeatletter\c@subsection=\uffa\makeatother

\subsection{Analytic functions over a $q$-invariant open disk}
\label{subsec:qtaylorconnesso}

In this section we prove proposition \ref{prop:taylor}.
Let us recall its statement:

\begin{prop}\label{prop:taylor2}
Let $D=D(\xi,\rho^-)$ be a $q$-invariant open disk.
The map
$$
\begin{array}{rccl}
T_{q,\xi}:&{\mathcal A}_D&\longrightarrow & K\{x-\xi\}_{q,\rho}\\ \\
          &f(x)          &\longmapsto\ds  &\sum_{n\geq 0}
               \ds \frac{d_q^nf}{[n]_q^!}(\xi)(x-\xi)_{q,n}
\end{array}
$$
is a $q$-difference algebras isomorphism.
Moreover, for all $f\in{\mathcal A}_D$, the series $T_{q,\xi}(f)(x)$
converges uniformly to $f(x)$ over any closed disk $D(\xi,\eta^+)$, with
$0<\eta<\rho$.
\end{prop}

\medskip
The proof is divided in three steps:

\medskip\noindent
{\it STEP 1. The map $T_{q,\xi}$ is a well-defined ring homorphism.}

\begin{proof}
The only non-trivial point is that $T_{q,\xi}$ is well-defined.
Let $f=\sum_{n\geq 0}f_n(x-\xi)^n\in{\mathcal A}_D$.
By (\ref{eq:normaqder}), for any integer $k\geq 0$
and any $\eta\in\R$ such that $|(1-q)\xi|\leq\eta<\rho$,
we have
\begin{equation*}\label{eq:norma}
    \l|\frac{d_q^k(f)}{[k]_q^!}(\xi)\r|
    \leq\l\|\frac{d_q^k}{[k]_q^!}f(x)\r\|_\xi(\eta)
    \leq\frac{\|f\|_\xi(\eta)}{\eta^k}\ .
\end{equation*}
Hence
$\liminf_{n\rightarrow\infty}\l|\frac{d_q^k(f)}{[k]_q^!}
(\xi)\r|^{-1/n}\geq \eta$, for all $\eta$ such that $|(1-q)\xi|\leq
\eta<\rho$.
Letting $\eta\rightarrow\rho$ one proves that $T_{q,\xi}(f)\in K\{x-\xi\}_{q,\rho}$.
\end{proof}

\medskip\noindent
{\it STEP 2. Let $y(x)=\sum_{n\geq 0}a_n(x-\xi)_{q,n}\in K\{x-\xi\}_{q,\rho}$.
Then $y(x)$ converges uniformly to an analytic function
over any closed disk $D(\xi,\eta^+)$, with
$0<\eta<\rho$.}

\begin{proof}
For any pair of integers $i$, $k$, such that $0\leq i\leq k$,
we consider the symmetric polynomial
$$
S^i_k(x_1,\dots,x_k)=
\begin{cases}
  1 &\hbox{if $i=0$;}\\
  \ds\sum_{1\leq j_1<\cdots<j_i\leq k} x_{i_1}\cdots x_{j_i}&\hbox{otherwise.}
\end{cases}
$$
We  notice that
\begin{equation}\label{eq:polsimmetrici}
\begin{array}{rcl}
(x-\xi)_{q,n}
   &=&(x-\xi)(x-q\xi)\cdots(x-q^{n-1}\xi)\\ \\
   &=&(x-\xi)\l[(x-\xi)+\xi(1-q)\r]\cdots\l[(x-\xi)+\xi(1-q^{n-1})\r]\\\\
   &=&\ds\sum_{k=0}^n\xi^{n-k}
      S^{n-k}_{n}(0, 1-q, \dots, 1-q^{n-1})(x-\xi)^k\ .
\end{array}
\end{equation}
Since
$|\xi^{n-k} S^{n-k}_{n}
(1-q^0, 1-q, \dots, 1-q^{n-1})|
\leq |\xi(1-q)|^{n-k-1}\leq R^{n-k-1}<\rho^{n-k-1}$,
the following series converges in $K$:
$$
\a_k=\sum_{n=k}^\infty a_n\xi^{n-k}
S^{n-k}_{n}(1-q^0, 1-q, \dots, 1-q^{n-1})\in K\ .
$$
Then the series $\widetilde y(x)\sum_{k\geq 0}\a_k(x-\xi)^k
\in K\[[x-\xi\]]$ converges over $D$, in fact for any $x_0\in D$ we have
$$
\sup\l(|\xi(1-q)|,|x_0-\xi|\r)\leq R<\rho\ ,
$$
and therefore
$$
\begin{array}{rcl}
0&\leq&\ds\limsup_{k\rightarrow\infty}|\a_k(x_0-\xi)^k|\\
 &\leq&\ds\limsup_{k\rightarrow\infty}\l|\sum_{n=k}^\infty a_n\xi^{n-k}
        S^{n-k}_{n}(1-q^0, 1-q,
        \dots, 1-q^{n-1})(x_0-\xi)^k\r|\\
 &\leq&\ds\limsup_{k\rightarrow\infty}\l(\sup_{n\geq k}
        \l|a_n\r||\xi(1-q)|^{n-k}|x_0-\xi|^k\r)\\
 &\leq&\ds\limsup_{k\rightarrow\infty}\l(\sup_{n\geq k}
        \l|a_n\r|R^n\r)=0
\end{array}\ .
$$
This proves that $\widetilde y(x)\in {\mathcal A}_D$.
\par
We want to prove that
$y(x)=\sum_{n\geq 0}a_n(x-\xi)_{q,n}\in K\{x-\xi\}_{q,\rho}$
converges uniformly to $\widetilde y(x)$ over any closed
disk $D(\xi,\eta^+)$, with $0<\eta<\rho$.
We fix a positive real number $\eta<\rho$
and we set $C_N=\sup_{n>N}|a_n|\eta^n$.
We observe that
$$
\lim_{N\rightarrow\infty} C_N=0\ ,
$$
since $\eta<\rho$.
We have
$$
\hskip -5pt
\begin{array}{rcl}
\lefteqn{\ds\Bigg\|\sum_{k=0}^N\a_k(x-\xi)^k-\sum_{n=0}^Na_k(x-\xi)_{q,k}\Bigg\|_\xi(\eta)}\\
&=&\ds\Bigg\|\sum_{k=0}^N\a_k(x-\xi)^k-\sum_{n=0}^N a_n\sum_{k=0}^n\xi^{n-k}
        S^{n-k}_{n}(0, 1-q, \dots, 1-q^{n-1})(x-\xi)^k
        \Bigg\|_\xi(\eta)\\
&=&\ds\Bigg\|\sum_{k=0}^N\sum_{n=k}^\infty a_n\xi^{n-k}
        S^{n-k}_{n}(0,1-q, \dots, 1-
        q^{n-1})(x-\xi)^k\\
&&\ds-\sum_{k=0}^N\sum_{n=k}^N a_n\xi^{n-k}
      S^{n-k}_{n}(0, 1-q, \dots, 1-q^{n-1})(x-\xi)^k\Bigg\|_\xi(\eta)\\
&=&\ds\Bigg\|\sum_{k=0}^N\Bigg(\sum_{n>N}a_n\xi^{n-k}
        S^{n-k}_{n}(1-q^0, 1-q, \dots, 1-q^{n-1})\Bigg)
        (x-\xi)^k\Bigg\|_\xi(\eta)\\
&\leq&\ds\sup_{k=0,\dots,N}\l(\eta^k\sup_{n>N}|a_n|\eta^{n-k}\r)\\
&\leq&\ds\sup_{k=0,\dots,N}\l(\sup_{n>N}|a_n|\eta^n\r)=C_N\ .
\end{array}
$$
We conclude that the two series converge to the same sum
over $D(\xi,\eta^+)$, for all $0<\eta<\rho$, hence over $D$.
\end{proof}

\medskip\noindent
{\it STEP 3. End of the proof.}

\begin{proof}
By the previous step, one can define a map
$S:K\{x-\xi\}_{q,\rho}\longrightarrow {\mathcal A}_D$.
Since $\xi$ is a limit point of $q^\N\xi$,
for any $f(x),g(x)\in{\mathcal A}_D$ we have
$$
\begin{array}{rcl}
\ds\frac{d_q^nf}{[n]_q^!}(\xi)
=\frac{d_q^ng}{[n]_q^!}(\xi)\ ,\ \forall n\in\Z_{n\geq 0}
&\Leftrightarrow&\hbox{$f(q^n\xi)=g(q^n\xi)$ for all integers $n\geq 0$}\\
&\Leftrightarrow&\hbox{$f=g$}\ .
\end{array}
$$
Then one deduces that $T_{q,\xi}^{-1}=S$ from the fact that
$$
\frac{d_q^n}{[n]_q^!}\l(S\circ T_{q,\xi}(f)\r)(\xi)
=\frac{d_q^nf}{[n]_q^!}(\xi)
\hbox{ and }
\frac{d_q^n}{[n]_q^!}\l(T_{q,\xi}\circ S (g)\r)(\xi)
=\frac{d_q^ng}{[n]_q^!}(\xi)
$$
for any $f\in{\mathcal A}_D$, any $g\in K\{x-\xi\}_{q,\rho}$
and any nonnegative integer $n$.
\end{proof}

This completes the proof of (\ref{prop:taylor2}).


\subsection{Analytic functions over non connected analytic domain}
\label{subsec:taylor}

Proposition \ref{prop:taylor2}
can be generalized to the
case of analytic functions over convenient
$q$-invariant non-connected analytic domains:

\begin{defn}
We call \emph{$q$-disk of center $\xi\in{\mathbb A}^1_K$ and radius
$\eta\in\R_{>0}$} the set
$$
q^\Z D(\xi,\eta^-):=\mathop\cup_{n\in\Z}D(q^n\xi,\eta^-)\ .
$$
\end{defn}

\begin{rmk}
Of course, if $n_0$ is the smallest positive integer such that
$|(1-q^{n_0})\xi|<\eta$, the $q$-disk $D$ is a disjoint union
of the open disks $D(q^n\xi,\eta^-)$, with $0\leq n<n_0$, each one of them
being $q^{n_0}$-invariant.
\end{rmk}

The algebra ${\mathcal A}_D$ of analytic functions over a $q$-disk $D$ is the
direct product of the algebras of analytic functions over each connected
component. Since the analytic domain $D$ is $q$-invariant ${\mathcal A}_D$ has a
structure of $q$-difference algebra.
Hence one can define a $q$-expansion map
$$
\begin{array}{rccl}
T_{q,\xi} & :{\mathcal A}_D & \longrightarrow &
    K\[[x-\xi\]]_q =\l\{\sum_{n\geq 0}a_n(x-\xi)_{q,n}~:~a_n\in K\r\}\\ \\
& f(x) & \longmapsto & \ds\sum_{n\geq 0}\frac{\dq^n(f)}{[n]_q^!}(\xi)(x-\xi)_{n,q}
\end{array}\ ,
$$
that we still call $T_{q,\xi}$.

\begin{prop}\label{prop:qdisk}
Let $f(x)=\sum_{n\geq 0}a_n(x-\xi)_{q,n}$ be a formal series such that the $a_n$'s
are elements of $K$
and $\ds\liminf_{n\rightarrow\infty}|a_n|^{-1/n}=\rho$.
Then $f(x)$ is a $q$-expansion of an analytic function over a $q$-disk
$D$ if and only if
$$
\frac{\rho}{|(q-1)\xi\pi_q|}>1\ .
$$
Let $n_0$ be the number of connected components of $D$ and $\widetilde\rho$ their radius.
Then $n_0$ is the smallest positive integer such that
\begin{equation}\label{eq:n_0}
\frac{\rho}{|(1-q)\xi\pi_q|}>|\pi_{q^{n_0}}|^{-1/n_0}
\end{equation}
and $\rho$, $\widetilde\rho$ and $n_0$ are linked by the relation:
\begin{equation}\label{eq:rhotilderho}
\frac{\widetilde\rho}{|(q^{n_0}-1)\xi\pi_{q^{n_0}}|}
=\l(\frac{\rho}{|(q-1)\xi\pi_q|}\r)^{n_0}\ .
\end{equation}
\end{prop}

\begin{rmk}
We recall that $\pi_q\in K $ satisfies
$|\pi_q|=\lim_{n\rightarrow\infty}|[n]_q^!|^{1/n}$
(cf. (\ref{defn:pi})).
\end{rmk}

It follows immediately by (\ref{cor:solanalytic}) and (\ref{prop:qdisk}) that:

\begin{cor}
The system
$$
Y(qx)=A(x)Y(x)\ ,
\hbox{ with $A(x)\in Gl_\mu({\mathcal M}_{D(\xi,\rho^-)})$,}
$$
has an analytic fundamental solution over a $q$-disk centered at $\xi$
if and only if the matrix $A(x)$ does not have any poles in $q^\N\xi$,
$\det A(x)$ does not have any zeros in $q^\N\xi$
and
$$
\limsup_{n\rightarrow\infty} \l|\frac{G_n(\xi)}{[n]_q^!}\r|^{1/n}<|(q-1)\xi\pi_q|^{-1}\ .
$$
\end{cor}

\subsection{Proof of proposition \ref{prop:qdisk}}
\label{subsec:qtaylor}

\medskip
The main result of this section is the proof proposition
\ref{prop:qdisk} (cf. (\ref{prop:provaqdisk}) below),
that we will deduce from the more general result (\ref{prop:polidisco}).

\medskip
Let $\xi\neq 0$ and $f(x)=\sum_{n\geq 0}a_n(x-\xi)_{q,n}\in K\[[x-\xi\]]_q$.
For any non-negative integer $k$ it makes sense to define formally
$$
 \frac{\dq^k(f)}{[k]_q^!}(\xi)=a_k
$$
and
$$
f(q^k\xi)=\sum_{n\geq 0}\a_n\xi^n(q^k-1)_{q,n}\ .
$$
An element $f(x)\in K\[[x-\xi\]]_q$ is uniquely determined
by $(f(q^k\xi))_{k\geq 0}$ or by $(a_k)_{k\geq 0}$, knowing that these two sequences are linked
by relations that can be deduced by (\ref{lemma:leibniz}).
Therefore $K\[[x-\xi\]]_q$ is not a local ring and we have:

\begin{lemma}
If $\xi\neq 0$ the natural morphism of $q$-difference algebras, defined for any $n_0\geq 1$,
\begin{equation}\label{eq:isoformale}
 \begin{array}{ccc}
 K\[[x-\xi\]]_q&\longrightarrow&\ds
       \prod_{i=0}^{n_0-1}K\[[x-q^i\xi\]]_{q^{n_0}}\\
 f(x)&\longmapsto &\ds
       \l(\sum_{n\geq 0}\frac{d_{q^{n_0}}^n(f)(q^i\xi)}{[n]_{q^{n_0}}^!}
       (x-q^i\xi)_{q^{n_0},n}\r)_{i=0,\dots,n_0-1}
 \end{array}
\end{equation}
is an isomorphism.
In particular  $K\[[x-\xi\]]_q$ is not a domain.
 \end{lemma}

\begin{proof}
It is enough to observe that an element $f\in\prod_{i=0}^{n_0-1}K\[[x-q^i\xi\]]_{q^{n_0}}$
is uniquely determined by
$$
\l(\l(f(q^{nn_0+i}\xi)\r)_{n\geq 0}\r)_{i=0,\dots,n_0-1}\ .
$$
\end{proof}

The isomorphism (\ref{eq:isoformale})
induces an isomorphism between ``converging $q$-series'':

\begin{prop}\label{prop:polidisco}
We set $\widetilde q=q^{n_0}$. Let
$$
g=\sum_{n\geq 0}g_n(x-\xi)_{q,n}\in K\[[x-\xi\]]_q
$$
and let
$$
f=\l(\sum_{n\geq 0}f_{i,n}(x-\xi)_{\tilde q,n}\r)_{0\leq i<n_0}
\in\prod_{i=0}^{n_0-1}K\[[x-q^i\xi\]]_{\tilde q}
$$
be its image via (\ref{eq:isoformale}).
We set
$$
\rho=\liminf_{n\rightarrow\infty}|g_n|^{-1}
\hbox{ and }
\ \widetilde\rho=\inf_{0\leq i<n_0}
\l(\liminf_{n\rightarrow\infty}|f_{i,n}|^{-1}\r)\ .
$$
Then
$$
\frac{\widetilde\rho}{|(\tilde q-1)\xi\pi_{\tilde q}|}
=\l(\frac{\rho}{|(q-1)\xi\pi_q|}\r)^{n_0}\ .
$$
\end{prop}

\medskip
In order to prove proposition
\ref{prop:polidisco}, we need a technical lemma:

\begin{lemma}
\par\noindent
1) Let $g\in K\[[x-\xi\]]_q$; then
\begin{equation}\label{eq:reluno}
  \begin{array}{rcl}
         \ds{(\widetilde q-1)^k\xi^k} d_{\widetilde q}^k(g)(\xi)
               &=&\ds (-1)^k\sum_{h=0}^{n_0k}
               \sum_{\frac{h}{n_0}\leq j\leq k}(-1)^j{k\choose j}_{\tilde q^{-1}}\\ \\
         &&\ds\cdot\widetilde q^{\frac{-j(j-1)}{2}}{n_0j\choose h}_q
               \cdot q^{\frac{h(h-1)}{2}} (q-1)^h\xi^h\dq^h(g)(\xi)\ .
  \end{array}
\end{equation}
\par\noindent
2) Let $f=(f_i)\in \prod_{i=0}^{n_0-1}K\[[x-q^i\xi\]]_{\tilde q}$. Then
\begin{equation}\label{eq:reldue}
  \begin{array}{rcl}
         \ds {(q-1)^k\xi^k}\dq^k(f)(\xi)
               &=&\ds(-1)^k\sum_{h=0}^{\l[\frac{n_0}{k}\r]}
               \sum_{\begin{array}{c}\scriptstyle 0\leq i<n_0
                 \\ \scriptstyle hn_0+i\leq k\end{array}}
               {k\choose hn_0+i}_{q^{-1}}q^{\frac{-(hn_0+i)(hn_0+i-1)}{2}}\\ \\
         &&\ds\cdot\sum_{j=0}^h{h\choose j}_{\tilde q}\widetilde q^{\frac{j(j-1)}{2}}q^{ij}
               (\widetilde q-1)^j\xi^jd_{\widetilde q}^j(f_i)(q^i\xi)\ .
  \end{array}
\end{equation}
\end{lemma}

\begin{proof}
The two formulas can be obtained by applying twice
the formulas relating $\dq$ and $\sgq$
(cf. \hbox{(\ref{lemma:leibniz})}),
taking into account that
$d_{\tilde q}(f(q^i\xi))=q^id_{\tilde q}(f)(q^i\xi)$.
\end{proof}

\begin{parag}
{\bf Proof of proposition \ref{prop:polidisco}.}

\par\medskip\noindent
The proof is based on the following property
of the limit of a sequence $(a_n)$ of real non-negative numbers
(cf. \cite[II, 1.8 and 1.9]{AB}):
\begin{equation}\label{eq:ab}
  \sup(1,\limsup_{n\rightarrow\infty}a_n^{1/n})=
  \limsup_{n\rightarrow\infty}\l(\sup_{s\leq n}a_s\r)^{1/n}\ .
\end{equation}

\smallskip\noindent
{\it STEP 1. Proof of (\ref{prop:polidisco}) when one of the (equivalent) conditions\\
- $\rho<|(q-1)\xi\pi_q|$,\\
- $\widetilde\rho<|(\widetilde q-1)\xi\pi_{\tilde q}|$\\
is satisfied.}
\par\noindent
Suppose first that $\rho<|(q-1)\xi\pi_q|$.
Since $\rho^{-1}|\pi_q|>1$, equality (\ref{eq:ab}) implies that
\begin{equation*}
  \begin{array}{rcl}
     \rho^{-1}|\pi_q|
     &=&\ds|\pi_q|\limsup_{n\rightarrow\infty}
                 \l|\frac{d_q^n(g)}{[n]_q^!}(\xi)\r|^{1/n}\\
         &\leq&\ds\limsup_{n\rightarrow\infty}
                 \l|d_q^n(g)(\xi)\r|^{1/n}\\
         &=&\ds\limsup_{n\rightarrow\infty}\l(\sup_{s\leq n} \l|\dq^s(g)(\xi)\r|\r)^{1/n}\ .
  \end{array}
\end{equation*}
Since $|(q-1)\xi\pi_q|\rho^{-1}>1$ we deduce from (\ref{eq:reldue}) that
\begin{equation*}
  \begin{array}{rcl}
         |(q-1)\xi\pi_q|\rho^{-1}
         &\leq&\ds\limsup_{n\rightarrow\infty}\l|(q-1)^n\xi^n\dq^n(g)(\xi)\r|^{1/n}\\
         &\leq&\ds\limsup_{n\rightarrow\infty}\l(\sup_{s\leq \l[\frac{n}{n_0}\r], l=0,\dots,n_0-1}
         \l|(\widetilde q-1)^s\xi^sd_{\widetilde q}^s(g)(q^l\xi)\r|\r)^{1/n}\\
         &=&\ds \sup\l(1,|(\widetilde q-1)\xi\pi_{\tilde q}|^{1/n_0}\widetilde\rho^{-1/n_0}\r)\\
         &=&\ds\l(|(\widetilde q-1)\xi\pi_{\tilde q}|\widetilde\rho\r)^{-1/n_0}
  \end{array}
\end{equation*}
and hence that
$\rho\geq |(q-1)\xi\pi_q|
\l(\frac{\tilde\rho}{|(\tilde q-1)\xi\pi_{\tilde q}|}\r)^{1/n_0}$.
Moreover the last inequality implies that
$\widetilde\rho<|(\widetilde q-1)\xi\pi_{\tilde q}|$.
\par
Let us suppose now that
$\widetilde\rho<|(\widetilde q-1)\xi\pi_{\tilde q}|$. Then
$$
\widetilde\rho^{-1}|\pi_{\tilde q}|
\leq\sup_{i=0,\dots,n_0-1}\l(\limsup_{n\rightarrow\infty}
\l(\sup_{s\leq n} \l|d_{\tilde q}^s(f)(q^i\xi)\r|\r)^{1/n}\r)\ .
$$
We deduce by (\ref{eq:reluno}) that
\begin{equation*}
  \begin{array}{rcl}
         \ds|(\widetilde q-1)\xi\pi_{\tilde q}|\widetilde\rho^{-1}
         &\leq&\ds\limsup_{n\rightarrow\infty}\l|d_{\tilde q}^n(f)(\xi)\r|^{1/n}\\
         &\leq&\ds\limsup_{n\rightarrow\infty}\l(\sup_{s\leq nn_0}
         \l|(q-1)^s\xi^s\dq^s(f)(\xi)\r|\r)^{1/n}\\  \\
         &\leq&\ds |(q-1)\xi|^{n_0}|p|^{n_0/(p-1)}\rho^{-n_0}\\
  \end{array}
\end{equation*}
and hence that
$\widetilde\rho\geq
\l(\frac{|(q-1)\xi\pi_q|}{|(\tilde q-1)\xi\pi_{\tilde q}|^{1/n_0}}\r)^{-n_0}\rho^{n_0}$.
Notice that this inequality implies that $\rho< |(q-1)\xi\pi_q|$

\smallskip\noindent
{\it STEP 2. Proof of (\ref{prop:polidisco}) when one of the (equivalent) conditions\\
- $\rho\geq|(q-1)\xi\pi_q|$,\\
- $\widetilde\rho\geq|(\widetilde q-1)\xi\pi_{\tilde q}|$\\
is satisfied.}
\par\noindent
We choose $\ga\in K$ such that
$$
|\ga|>|(q-1)\xi\pi_q|^{-1}\rho\ .
$$
By setting $\xi=\ga\zeta$ and $x=\ga t$, we identify $g$ with
an element of $K\[[t-\zeta\]]_q$ satisfying the hypothesis of Step 1.
We deduce from the estimate in Step 1 that
\begin{equation*}
  \begin{array}{rcl}
         \ds|\ga|^{-1}\rho
         &=&\ds\frac{|(q-1)\zeta\pi_q|}{|(\widetilde q-1)\zeta\pi_{\tilde q}|^{1/n_0}}
       (|\ga|^{-1}\widetilde \rho)^{1/n_0}\\
         &=&\ds|\ga|^{-1} \frac{|(q-1)\xi\pi_q|}{|(\widetilde q-1)\xi\pi_{\tilde q}|^{1/n_0}}
           \widetilde\rho^{1/n_0}\ .
  \end{array}
\end{equation*}
\hfill$\square$
\end{parag}

\begin{parag}\label{prop:provaqdisk}
{\bf Proof of proposition \ref{prop:qdisk}.}

\par\medskip\noindent
We briefly recall the statement of (\ref{prop:qdisk}).
Let $f(x)=\sum_{n\geq 0}a_n(x-\xi)_{q,n}$ be a formal series such that $a_n\in K$
and $\ds\liminf_{n\rightarrow\infty}|a_n|^{-1/n}=\rho$.
We have to prove that $f(x)$ is a $q$-expansion of an analytic function over a $q$-disk
$q^\Z D(\xi,\widetilde\rho^-)$ if and only if
\begin{equation}\label{eq:hyp}
\frac{\rho}{|(q-1)\xi\pi_q|}>1\ .
\end{equation}
Let $f(x)$ satisfy (\ref{eq:hyp}).
Since $\frac{\rho}{|(q-1)\xi\pi_q|}>1$ there exists a minimal
positive integer $n_0$ such that
$$
\l(\frac{\rho}{|(q-1)\xi\pi_q|}\r)^{n_0}|\pi_{q^{n_0}}|>1\ .
$$
Then set
$$
\widetilde\rho
=\l(\frac{\rho}{|(q-1)\xi\pi_q|}\r)^{n_0}|\pi_{q^{n_0}}|{|(q^{n_0}-1)\xi|}\
.
$$
It follows from (\ref{prop:taylor}) and (\ref{prop:qdisk})
that $f(x)$ is the $q$-expansion
of an analytic function over $q$-disk $q^\Z D(\xi,\widetilde\rho^-)$
satisfying (\ref{eq:n_0}) and (\ref{eq:rhotilderho}).
\par
Conversely, let $f(x)$ be the $q$-expansion of an analytic function over a $q$-disk
$D=q^\Z D(\xi,\widetilde\rho^-)$.
Let $n_0$ be the number of connected components of $D$
and $\widetilde q=q^{n_0}$. Then $\widetilde\rho>|(\widetilde q-1)\xi|$.
It follows from (\ref{prop:polidisco}) that $\rho>|(q-1)\xi\pi_q|$.
\hfill$\square$
\end{parag}

\makeatletter\newcount\uffa\uffa=\c@subsection\makeatother
\section{Basic facts about regular singularities of $q$-dif\-fer\-ence systems}
\label{sec:regsing}
\makeatletter\c@subsection=\uffa\makeatother

In this chapter we briefly recall some properties of regular singular
$q$-difference systems, that we used in \ref{sec:frobenius} and
\ref{sec:transfer}. A complete exposition can be found in
\cite[Ch. 12]{PS}, \cite{Sfourier} and \cite{Spoligono}.

\subsection{Regular singular $q$-difference systems}

Let $K$ be an algebraically closed field of characteristic $0$ and
$K((x))=K\[[x\]]\l[\frac{1}{x}\r]$ the field of Laurent series.
For any $q\in K$, it has a natural structure of $q$-difference algebra,
hence one can consider a $q$-difference system with coefficients in
$K((x))$:
$$
Y(qx)=A(x)Y(x)\ .
$$

\begin{defn}
The $q$-difference system $Y(qx)=A(x)Y(x)$ is said to have a
\emph{regular singularity at $0$}
(or is \emph{regular singular at $0$}) if there exists
$U(x)\in Gl_\mu(K((x)))$ such that the matrix $A_{[U]}(x)=U(qx)A(x)U(x)^{-1}\in Gl_\mu(K\[[x\]])$.
\end{defn}

\begin{lemma}{\rm\cite[1.1.3]{Sfourier}}\label{lemma:uniforme}
Suppose $A(x)\in Gl_\mu(K\[[x\]])$.
If any two eigenvalues $\a,\be$ of $A(0)$ are such that either
$\a=\be$ or $\a\be^{-1}\not\in q^\Z$, then one can construct $U(x)\in Gl_\mu(K\[[x\]])$ such that
$A_{[U]}(x)=A(0)$.
\end{lemma}

By a convenient gauge transformation one can always assume that the hypothesis of the
lemma above are satisfied:

\begin{prop}{\rm\cite[1.1.1]{Sfourier}}\label{prop:shearing}
Suppose that $A(x)\in Gl_\mu(K\[[x\]])$ and for any eigenvalue $\a$ of $A(0)$ choose an integer $n_\a$.
Then there exists a matrix $H(x)$ constructed by alternatively multiplying constant matrices
in $Gl_\mu(K)$ and diagonal matrices of the form
\begin{equation}\label{eq:shearing}
\begin{pmatrix}
{\mathbb I}_{\mu_1}&                    &   \\
                   &{\mathbb I}_{\mu_2}x^{\pm 1}&   \\
                   &                    &{\mathbb I}_{\mu_3}\end{pmatrix}\ ,
\hbox{ with }\mu_1+\mu_2+\mu_3=\mu\ ,
\end{equation}
such that, for any eigenvalue $\a$ of $A(0)$,
$q^{n_\a}\a$ is an eigenvalue of $A_{[H]}(0)$.
\par
In particular one can construct a matrix $H(x)$ such that any two eigenvalues $\a,\be$ of
$A_{[H]}(0)$ satisfy either $\a=\be$ or $\a\be^{-1}\not\in q^\Z$.
\end{prop}

The idea of the proof of the proposition above is that one has first to
consider a constant gauge matrix $Q$ such that $Q^{-1}A(x)Q$ has the constant term in the
Jordan normal form. Then by using a gauge matrix of the form (\ref{eq:shearing}) one multiplies the eigenvalue of
a chosen block by $q^{\pm 1}$ of $A(0)$.
By iterating the algorithm one obtain the desired gauge transformation.
One calls a gauge
matrix constructed as in the previous proposition \emph{shearing transformation}.

\begin{cor}\label{cor:uniforme}
Let $U(x),V(x)\in Gl_\mu(K((x)))$ be two gauge transformation matrix such that
both $A_{[U]}$ and $A_{[V]}$ are in $GL_\mu(K)$. Then there exists a shearing transformation $H(x)$
such that  $HU=V$.
\end{cor}

The corollary follows immediately by proposition (\ref{prop:shearing}),
by observing that the Jordan normal forms of $A_{[U]}$ and $A_{[V]}$ coincide, modulo the fact that
the eigenvalues of $A_{[U]}$ are the eigenvalues of $A_{[V]}$ multiplied by an integer
power of $q$.

\subsection{From $q$-difference systems to $q$-difference equations}
\label{subsec:qequation}

Consider now a $q$-difference equation
$$
{\mathcal L}y=a_\mu(x)y(q^\mu x)-a_{\mu-1}(x)y(q^{\mu-1}x)-\dots-a_0(x)y(x)=0\ ,
$$
with $a_i(x)\in K((x))$, for all $i=0,\dots,\mu-1$, and $a_\mu(x)=1$.
The origin is said to be a \emph{regular singularity of ${\mathcal L}y=0$}
if and only if the \emph{Newton Polygon}, i.e. the convex hull in $\R^2$ of
$$
\{(i,j)\in\Z^2:\
i=0,\dots,\mu \hbox{ and } j\geq{\rm ord}_xa_i(x)\}\ ,
$$
has only one finite slope equal to $0$.
Obviously one has:

\begin{lemma}
If $0$ is a regular singularity for ${\mathcal L}y=0$, then
$a_i(x)$ does not have any pole at $0$, for any $i=0,\dots,\mu-1$.
In particular, $a_0(0)\neq 0$.
\end{lemma}

It follows that the $q$-difference system
\begin{equation}
Y(qx)=A_{[H]}(x)Y(x)\ ,
\hbox{ with }A_{[H]}(x)=
\l(\begin{array}{c|c}
0& \\
\vdots &{\mathbb I}_{\mu-1}\\
0& \\
\hline\\
a_0(x)&a_1(x)\dots a_{\mu-1}(x)
\end{array}\r)\
\end{equation}
has a regular singularity at $0$.

\medskip
The converse is also true:

\begin{prop}{\rm(\cite[Annexe B]{Sfourier} and \cite[\S2, in particular 2.2.6, (ii)]{Spoligono})}\\
A $q$-difference system over $K((x))$ has a regular singularity at $0$
if and only if the associated $q$-difference equation \emph{via} the cyclic vector lemma
has a regular singularity at $0$.
\end{prop}

\begin{parag}\label{parag:secorder}
{\bf Second order $q$-difference equations.}

\par\medskip\noindent
Consider a second order regular singular $q$-difference equation
\begin{equation}\label{eq:secorder}
y(q^2x)-P(x)y(qx)-Q(x)y(x)=0\ .
\end{equation}
Then the associated $q$-difference system is given by
$$
Y(qx)=A(x)Y(x)\ ,
\hbox{ with }
A(x)=\begin{pmatrix}0 & 1\\Q(x) &P(x)\end{pmatrix}\
$$
and the eigenvalues $\a,\be$ of $A(0)$ are solutions of the second order equation

\begin{equation}\label{eq:eqindiciale}
T^2-P(0)T-Q(0)=0\ .
\end{equation}
By (\ref{lemma:uniforme}), if $\a\neq\be$ and $\a\be^{-1}\not\in q^\Z$, there exists
$U(x)\in Gl_\mu(K\[[x\]])$ such that
$$
U(qx)A(x)U(x)^{-1}=\begin{pmatrix}\a&0\\0&\be\end{pmatrix}\ .
$$
If we call $e_\a(x)$ (resp. $e_\be(x)$) a solution of $y(qx)=\a y(x)$
(resp. $y(qx)=\be y(x)$) in a convenient $q$-difference extension of $K((x))$,
then we have
$$
\l[U(qx)^{-1}\begin{pmatrix}e_\a(qx)&0\\0&e_\be(qx)\end{pmatrix}\r]
=A(x)\l[U(x)^{-1}\begin{pmatrix}e_\a(x)&0\\0&e_\be(x)\end{pmatrix}\r]\ .
$$
If $(u_\a(x),u_\be(x))$ is the first row of $U(x)^{-1}$, then
$e_\a(x)u_\a(x),e_\be(x)u_\be(x)$ is a basis of solutions of (\ref{eq:secorder}),
meaning that $e_\a(x)u_\a(x)$ and $e_\be(x)u_\be(x)$ are
linearly independent over the field of constants and they span the
vector space of solutions of ${\mathcal L}y=0$.
Observe also that $e_\a(x),e_\be(x)$ must be linearly independent over $K((x))$
(that can be proved as in \cite[Annexe A, 4)]{Sfourier}).
\end{parag}

\makeatletter\newcount\uffa\uffa=\c@subsection\makeatother
\section{The $q$-type of a number}
\label{sec:type}
\makeatletter\c@subsection=\uffa\makeatother

The purpose of this appendix is to calculate the radius of convergence
of the series
$$
\Phi(x)=\sum_{n\geq 0}\frac{(1-q)^nx^n}{(1-q\a)\cdots(1-q^n\a)}
$$
for $\a\not\in q^{-\N}$, which is a key point in the proof theorem \ref{thm:transfer}.
Actually the radius of such a series turns out to be equal to $|\pi|{\rm type}_q(\a)$.

\medskip
As in chapters \ref{sec:frobenius} and \ref{sec:transfer},
we assume that $|1-q|<|\pi|$, so that $|1-q|=|\log q|$.

\subsection{Basic properties of the $q$-type of a number}
\label{subsec:type}

We recall that for any $\a\in K$ we have defined the $q$-type of $\a$ to be
the radius of convergence of $\sum_{n\geq 0\atop\a\neq q^n}\frac{(1-q)^nx^n}{1-q^n\a}$, while the
type of $\a$ is the radius of convergence of $\sum_{n\geq 0\atop\a\neq n}\frac{x^n}{n-\a}$.

\begin{prop}\label{prop:qtype}
For any $\a\in K$ we have
$$
{\rm type}_q(\a)=\l\{\begin{array}{ll}
\ds{\rm type}\l(\frac{\log\a}{\log q}\r)&\hbox{ if }\ds\l|\frac{\a-1}{q-1}\r|\leq 1\ ,\\ \\
1&\hbox{ otherwise.}
\end{array}\r.
$$
\end{prop}

\begin{proof}
We have
$$
\frac{\a-q^n}{q-1}=\frac{\a-1}{q-1}-\frac{q^n-1}{q-1}=\frac{\a-1}{q-1}-[n]_q\ .
$$
Hence, if $\l|\frac{\a-1}{q-1}\r|>1$, the radius of convergence of $\sum_{n\geq 0}\frac{q-1}{\a-q^n}x^n$
is equal to $1$.
\par
Suppose $|\a-1|\leq|q-1|<|\pi|$. Then $|\a-q^n|\leq|q-1|<|\pi|$ for any integer $n$ and
$$
\l|\frac{\a-q^n}{q-1}\r|
=\l|\frac{\log(\a q^{-n})}{\log q}\r|
=\l|\frac{\log\a}{\log q}-n\r|\ .
$$
\end{proof}

\subsection{Radius of convergence of ${}_1\Phi_1(q;\a q;q;(1-q)x)$}

The series $\Phi(x)$,
whose radius of convergence we want to estimate, is a basic hypergeometric series.
In the literature it is denoted by
${}_1\Phi_1(q;\a q;q;(1-q)x)$ (cf. \cite{GR}).

\begin{lemma}
For any $\a\in K\smallsetminus q^{\Z}$ we have:
$$
\Phi(x)=
\frac{1-\a}{1-q}\l(\sum_{n\geq 0}\frac{x^n}{[n]_q^!}\r)
\l(\sum_{n\geq 0}q^{\frac{n(n+1)}{2}}\frac{(-x)^n}{[n]_q^!}\frac{1-q}{1-q^n\a}\r)\ .
$$
\end{lemma}

\begin{proof}
One verifies directly that the series $\Phi(x)$
is a solution of the $q$-difference equation
$$
\begin{array}{rcl}
{\mathcal L}\Phi(x)
&=&\big[\sgq-1\big]\circ\big[\a\sgq-\l((q-1)x+1\r)\big]\Phi(x)\\ \\
&=&\a\Phi(q^2x)-\l((q-1)qx+1+\a\r)\Phi(qx)+(q-1)qx+1=0\ .
\end{array}
$$
Since the roots of the equation (cf. (\ref{eq:eqindiciale}))
$$
\a T^2-(\a+1)T+1=0
$$
are exactly $\a^{-1}$ and $1$, any solution of ${\mathcal L}y(x)=0$ of the
form $1+\sum_{n\geq 1}a_nx^n\in K\[[x\]]$ must coincide with $\Phi(x)$.
Therefore, to finish the proof of the lemma, it is enough to verify that
$$
\Psi(x)=\frac{1-\a}{1-q}\l(\sum_{n\geq 0}\frac{x^n}{[n]_q^!}\r)
\l(\sum_{n\geq 0}q^{\frac{n(n+1)}{2}}\frac{(-x)^n}{[n]_q^!}\frac{1-q}{1-q^n\a}\r)
$$
is a solution of ${\mathcal L}y(x)=0$ and that $\Psi(0)=1$.
\par
Let $e_q(x)=\sum_{n\geq 0}\frac{x^n}{[n]_q^!}$. Then $e_q(x)$ satisfies the $q$-difference equation
$$
e_q(qx)=\l((q-1)x+1\r)e_q(x)\ ,
$$
hence
$$
\begin{array}{rcl}{\mathcal L}\circ e_q(x)
&=&\big[\sgq-1\big]\circ e_q(qx)\circ \big[\a\sgq-1\big]\\ \\
&=&e_q(x)\l((q-1)x+1\r)\big[\l((q-1)qx+1\r)\sgq-1\big]\circ \big[\a\sgq-1\big]\\ \\
&=&(\ast)\big[\l((q-1)qx+1\r)\sgq-1\big]\circ \big[\a\sgq-1\big]\ ,
\end{array}
$$
where we have denoted with $(\ast)$ a coefficient in $K(x)$, not depending on $\sgq$.
\par
Consider the series $E_q(x)=\sum_{n\geq 0}q^{\frac{n(n-1)}{2}}\frac{x^n}{[n]_q^!}$, which satisfies
$$
\l(1-(q-1)x\r)E_q(qx)=E_q(x)\
$$
and the series
$$
g_\a(x)=\sum_{n\geq 0}q^{\frac{n(n+1)}{2}}\frac{(-x)^n}{[n]_q^!}\frac{1-q}{1-q^n\a}\ .
$$
Then
$$
\begin{array}{rcl}{\mathcal L}\circ e_q(x)g_\a(x)
&=&(\ast)\big[\l((q-1)qx+1\r)\sgq-1\big]\circ \big[\a\sgq-1\big]g_\a(x)\\ \\
&=&(\ast)\big[\l((q-1)qx+1\r)\sgq-1\big] E_q(-qx)\\ \\
&=&(\ast)\big[\l((q-1)qx+1\r)E_q(-q^2x)-E_q(-qx)\big]\\ \\
&=&0\ .
\end{array}
$$
It is enough to observe that $e_q(0)g_\a(0)=\frac{1-q}{1-\a}$
to conclude that the series $\Psi(x)=\frac{1-\a}{1-q}e_q(x)g_\a(x)$ coincides with $\Phi(x)$.
\end{proof}

\begin{cor}\label{cor:type}
For any $\a\in K\smallsetminus q^{\Z}$, the radius of convergence of
$$
\sum_{n\geq 0}\frac{(1-q)^nx^n}{(1-q\a)\cdots(1-q^n\a)}
$$
is $|\pi|{\rm type}_q(\a)$.
\end{cor}

\begin{proof}
It follows immediately from the previous results, since the radius of convergence of
$\sum_{n\geq 0}\frac{x^n}{[n]_q^!}$ is $|\pi|$.
\end{proof}

\end{appendix}

\bigskip\noindent
{\bf Acknowledgments.}
As an undergraduate student in Padova, I attended Professor Dwork's course
during the academic year 1993/94: it is an honor to be able to acknowledge
his decisive influence on my Ph.D. thesis as well as on the present work.
\par
I would like to thank
Y. Andr\'e, J-P. Ramis and J. Sauloy for the many discussions we have had
during the preparation of this paper, for sharing their enthusiasm, and for
their constant encouragement; the referee of the present paper for writing
such a motivating report; and H. Shapiro and E. Sommers for their editorial help.
Any remaining deficiencies are entirely my fault.
\frenchspacing

\newcommand{\noopsort}[1]{}

\makelastpage

\end{document}